\documentclass[12pt,twoside]{amsart}
\usepackage{amssymb,amscd,amsxtra,calc,mathrsfs}
\usepackage{amsmath}
\usepackage{xypic}

\title[Uniform K-stability of $\mathbb{Q}$-Fano varieties]{A valuative criterion for 
uniform K-stability of $\mathbb{Q}$-Fano varieties}
\author{Kento Fujita} 
\date{\today}
\address{Research Institute for Mathematical Sciences, Kyoto University, Kyoto 606-8502, Japan}
\email{fujita@kurims.kyoto-u.ac.jp}

\newcommand{\pr}{\mathbb{P}}

\newcommand{\Z}{\mathbb{Z}}
\newcommand{\Q}{\mathbb{Q}}
\newcommand{\R}{\mathbb{R}}
\newcommand{\C}{\mathbb{C}}

\newcommand{\A}{\mathbb{A}}
\newcommand{\G}{\mathbb{G}}

\newcommand{\Supp}{\operatorname{Supp}}
\newcommand{\Exc}{\operatorname{Exc}}

\newcommand{\DIV}{\operatorname{div}}

\newcommand{\lct}{\operatorname{lct}}
\newcommand{\DF}{\operatorname{DF}}
\newcommand{\Ding}{\operatorname{Ding}}
\newcommand{\ord}{\operatorname{ord}}

\newcommand{\vol}{\operatorname{vol}}
\newcommand{\Image}{\operatorname{Image}}

\newcommand{\CAN}{\operatorname{can}}
\newcommand{\LC}{\operatorname{lc}}
\newcommand{\AC}{\operatorname{ac}}
\newcommand{\NA}{\operatorname{NA}}
\newcommand{\MA}{\operatorname{MA}}
\newcommand{\s}{\operatorname{s}}
\newcommand{\red}{\operatorname{red}}

\newcommand{\triv}{\operatorname{triv}}
\newcommand{\sI}{\mathcal{I}}

\newcommand{\sO}{\mathcal{O}}

\newcommand{\sX}{\mathcal{X}}
\newcommand{\sY}{\mathcal{Y}}
\newcommand{\sL}{\mathcal{L}}

\newcommand{\sH}{\mathcal{H}}

\newcommand{\sF}{\mathcal{F}}

\newcommand{\sZ}{\mathcal{Z}}

\newtheorem{thm}{Theorem}[section]

\newtheorem{proposition}[thm]{Proposition}
\newtheorem{corollary}[thm]{Corollary}
\newtheorem{claim}[thm]{Claim}

\theoremstyle{definition}
\newtheorem{definition}[thm]{Definition}
\newtheorem{remark}[thm]{Remark}

\newtheorem{example}[thm]{Example}
\newtheorem*{ack}{Acknowledgments}

\begin{document}

\maketitle 

\begin{abstract}
We give a simple necessary and sufficient condition for uniform K-stability 
of $\mathbb{Q}$-Fano varieties. 
\end{abstract}

\setcounter{tocdepth}{1}
\tableofcontents

\section{Introduction}\label{intro_section}

Throughout the article, we work out over an algebraically closed field $\Bbbk$ of 
characteristic zero unless otherwise noted. 
A \emph{$\Q$-Fano variety} is defined by a normal projective 
variety $X$ which is log terminal and such that 
the anti-canonical divisor $-K_X$ is an ample 
$\Q$-Cartier divisor. A \emph{complex Fano manifold} is defined by a 
$\Q$-Fano variety over the complex number field $\C$ which is smooth over $\C$. 
Take any $\Q$-Fano variety $X$. It is important to know whether $X$ is 
\emph{K-stable}, or \emph{K-semistable}, or not. Indeed, if $X$ is a 
complex Fano manifold, then $X$ admits K\"ahler-Einstein metrics if and only if 
$X$ is \emph{K-polystable} (see \cite{tian1, don05, CT, stoppa, mab1, mab2, B, 
CDS1, CDS2, CDS3, tian2}). K-stability is stronger than K-polystability, and 
K-polystability is stronger than K-semistability. 
In this article, we mainly treat \emph{uniform K-stability} which is stronger than 
K-stability. The notion of uniform K-stability was originally introduced by Sz\'ekelyhidi 
\cite{sz1, sz2} and was deeply developed in \cite{dervan, BHJ}.
Moreover, if $X$ is a complex Fano manifold, then uniform K-stability of $X$ is 
equivalent to K-stability of $X$ 
by \cite{CDS1, CDS2, CDS3, tian2} and \cite{BBJ}. 
We remark that the holomorphic automorphism group of $X$ is finite if 
$X$ is uniformly K-stable; see \cite[Corollary E]{BHJ2} for details. 
We will define those stability notions in Section \ref{tc_section}. 

We show that we can test uniform K-stability and K-semistability of $X$
by calculating certain invariants for all divisorial valuations on $X$. 

\begin{definition}\label{intro_dfn1}
Let $X$ be a $\Q$-Fano variety of dimension $n$. Take any projective birational 
morphism $\sigma\colon Y\to X$ with $Y$ normal and any prime divisor $F$ on $Y$, 
that is, $F$ is a prime divisor \emph{over} $X$ (see \cite[Definition 2.24]{KoMo}).
Take any $k\in\Z_{\geq 0}$ with $-kK_X$ Cartier and take any $x\in\R_{\geq 0}$. 
\begin{enumerate}
\renewcommand{\theenumi}{\arabic{enumi}}
\renewcommand{\labelenumi}{(\theenumi)}
\item\label{intro_dfn11}
We define the sub vector space $H^0(X, -kK_X-xF)$ of the $\Bbbk$-vector space 
$H^0(X, -kK_X)$ as 
\begin{eqnarray*}
&&H^0(X, -kK_X-xF):=H^0(Y, \sigma^*\sO_X(-kK_X)(\lfloor -xF\rfloor))\\
&\subset&H^0(Y, \sigma^*\sO_X(-kK_X))=H^0(X, \sO_X(-kK_X)).
\end{eqnarray*}
\item\label{intro_dfn12}
We define 
\begin{eqnarray*}
&&\vol_X(-K_X-xF):=\vol_Y(\sigma^*(-K_X)-xF)\\
&=&\limsup_{\substack{k\to\infty\\-kK_X:\text{ Cartier}}}
\frac{\dim H^0(X, -kK_X-kxF)}{k^n/n!}.
\end{eqnarray*}
By \cite[Corollary 2.2.45]{L1} and \cite[Example 11.4.7]{L2}, 
the limsup is actually a limit and the function $\vol_X(-K_X-xF)$ 
is a monotonically decreasing continuous function for $x\in\R_{\geq 0}$. 
\end{enumerate}
\end{definition}

\begin{remark}\label{intro_rmk}
The above invariants do not depend on the choice of $\sigma\colon Y\to X$. 
More precisely, these invariants 
depend only on the divisorial valuation on $X$ given by $F$. 
See also Section \ref{val_section}. 
\end{remark}

\begin{definition}\label{intro_dfn2}
Let $X$ be a $\Q$-Fano variety and $F$ be a prime divisor over $X$. 
\begin{enumerate}
\renewcommand{\theenumi}{\arabic{enumi}}
\renewcommand{\labelenumi}{(\theenumi)}
\item\label{intro_dfn20}
The divisor $F$ is said to be \emph{dreamy}\footnote{The word ``dreamy" comes 
from ``Mori dream spaces" in \cite{HK}.} if the $\Z_{\geq 0}^{\oplus 2}$-graded 
$\Bbbk$-algebra 
\[
\bigoplus_{k,j\in\Z_{\geq 0}}H^0(X, -krK_X-jF)
\] 
is finitely generated for some $r\in\Z_{>0}$ with $-rK_X$ Cartier. 
(The definition does not depend on the choice of $r$.) 
\item\label{intro_dfn21}
We define the \emph{pseudo-effective threshold} $\tau(F)\in\R_{>0}$ \emph{of} 
$F$ \emph{with respect to} $-K_X$ as 
\[
\tau(F):=\sup\{x\in\R_{>0}\,|\, \vol_X(-K_X-xF)>0\}.
\]
\item\label{intro_dfn22}
We define the \emph{log discrepancy} $A_X(F)\in\Q_{>0}$ \emph{of} $F$ 
\emph{with respect to} $X$ as 
$A_X(F):=\ord_F(K_{Y/X})+1$. Since $X$ is log terminal, $A_X(F)>0$ holds.
\item\label{intro_dfn23}
We set 
\[
\beta(F):=A_X(F)\cdot\vol_X(-K_X)-\int_0^{\tau(F)}\vol_X(-K_X-xF)dx.
\]
\item\label{intro_dfn24}
We set 
\[
j(F):=\int_0^{\tau(F)}(\vol_X(-K_X)-\vol_X(-K_X-xF))dx.
\]
Since $\vol_X(-K_X)-\vol_X(-K_X-xF)>0$ for any $x\in(0, \tau(F))$ 
(see \cite[Lemma 5.13]{BHJ}), we have $j(F)>0$.
\end{enumerate}
\end{definition}

The following is the main theorem in this article (see also Theorem \ref{log_mainthm}).

\begin{thm}[Main Theorem]\label{mainthm}
Let $X$ be a $\Q$-Fano variety. For any $\delta\in[0, 1)$, the following are 
equivalent: 
\begin{enumerate}
\renewcommand{\theenumi}{\roman{enumi}}
\renewcommand{\labelenumi}{(\theenumi)}
\item\label{mainthm_1}
$\DF(\sX, \sL)\geq \delta\cdot J^{\NA}(\sX, \sL)$ holds for any normal, ample 
test configuration $(\sX, \sL)/\A^1$ for $(X, -K_X)$ $($see \S \ref{tc_section} 
for the definitions$)$.
\item\label{mainthm_2}
$\beta(F)\geq \delta\cdot j(F)$ holds for any prime divisor $F$ over $X$.
\item\label{mainthm_3}
$\beta(F)\geq \delta\cdot j(F)$ holds for any dreamy prime divisor $F$ over $X$.
\end{enumerate}
\end{thm}

As a direct consequence of Theorem \ref{mainthm}, we get the following corollary. 

\begin{corollary}\label{maincor}
Let $X$ be a $\Q$-Fano variety. 
\begin{enumerate}
\renewcommand{\theenumi}{\arabic{enumi}}
\renewcommand{\labelenumi}{(\theenumi)}
\item\label{maincor1}
The following are equivalent: 
\begin{enumerate}
\renewcommand{\theenumii}{\roman{enumii}}
\renewcommand{\labelenumii}{(\theenumii)}
\item\label{maincor11}
$X$ is uniformly K-stable, 
\item\label{maincor12}
there exists $\delta\in(0, 1)$ such that 
$\beta(F)\geq \delta\cdot j(F)$ holds for any prime divisor $F$ over $X$,
\item\label{maincor13}
there exists $\delta\in(0, 1)$ such that 
$\beta(F)\geq \delta\cdot j(F)$ holds for any \emph{dreamy} prime divisor $F$ over $X$.
\end{enumerate}
\item\label{maincor2}
The following are equivalent: 
\begin{enumerate}
\renewcommand{\theenumii}{\roman{enumii}}
\renewcommand{\labelenumii}{(\theenumii)}
\item\label{maincor21}
$X$ is K-semistable, 
\item\label{maincor22}
$\beta(F)\geq 0$ holds 
for any prime divisor $F$ over $X$,
\item\label{maincor23}
$\beta(F)\geq 0$ holds 
for any \emph{dreamy} prime divisor $F$ over $X$.
\end{enumerate}
\end{enumerate}
\end{corollary}

Corollary \ref{maincor} is important and useful in order to 
test K-semistability and uniform K-stability of $\Q$-Fano varieties. 
For example, we can very easily check 
that $\pr^1$ is K-semistable from the criterion (see also Example 
\ref{logP1_ex}). Moreover, from Corollary 
\ref{maincor} \eqref{maincor2}, we can immediately get the result in 
\cite[Theorem 1.1]{fjt2} (see the proof of \cite[Theorem 5.1]{fjt2}). 
Furthermore, in \cite{FO}, the author and Yuji Odaka recovered the proofs of 
various known criteria for 
K-semistability and uniform K-stability of $\Q$-Fano varieties 
by using Corollary \ref{maincor}. 

We also get a criterion for K-stability of $\Q$-Fano varieties. 

\begin{thm}\label{mainK_thm}
Let $X$ be a $\Q$-Fano variety. Then the following are equivalent: 
\begin{enumerate}
\renewcommand{\theenumi}{\roman{enumi}}
\renewcommand{\labelenumi}{(\theenumi)}
\item\label{mainK_thm1}
$X$ is K-stable, 
\item\label{mainK_thm2}
$\beta(F)>0$ holds for any \emph{dreamy} prime divisor $F$ over $X$. 
\end{enumerate}
\end{thm}

In the theory of birational geometry, it is important to consider various invariants 
of prime divisors over $X$. The criteria in Theorems \ref{mainthm} and \ref{mainK_thm} 
seem very natural and simple from the view point of birational geometry.

\begin{remark}\label{li_rmk}
\begin{enumerate}
\renewcommand{\theenumi}{\arabic{enumi}}
\renewcommand{\labelenumi}{(\theenumi)}
\item\label{li_rmk1}
Demanding that $\beta(F)>0$ (resp.\ $\beta(F)\geq 0$) for all prime divisors $F$ 
\emph{on} $X$ leads to be the (weaker) notion of divisorial stability 
(resp.\ divisorial semistability) studied in \cite{fjt1}. 
Note that every prime divisor on $X$ is dreamy by \cite[Corollary 1.3.2]{BCHM}.  
\item\label{li_rmk2}
It seems difficult to show $\beta(F)>0$ for \emph{any} prime divisor $F$ over $X$ 
under the assumption that $X$ is K-stable. Indeed, we consider a kind of the 
limit of certain invariants in order to show the positivity of $\beta(F)$ 
(see Section \ref{Kv_section}).
\item\label{li_rmk3}
If $x\in X$ is a general point and $F$ is the exceptional divisor of the 
blowup along $x\in X$, then the invariant $\beta(F)$ is closely related to the 
invariant $\beta_x(-K_X)$ in \cite[\S 4]{MR} which was first studied by 
Per Salberger in his unpublished work. Moreover, the invariant $j(F)$ essentially 
appeared in \cite[p.\ 407]{sal}. 
\item\label{li_rmk4}
Recently, Chi Li showed in \cite{li2} that K-semistability of $X$ 
is equivalent to the 
condition that the normalized volume of $\G_m$-invariant valuations 
over the affine cone $(C_{r_0}, 0)$ of $X$ with respect to some positive multiple 
$-r_0K_X$ of $-K_X$ is minimized at the canonical valuation. See also \cite{li1}. 
Corollary \ref{maincor} \eqref{maincor2} looks similar to \cite[Theorem 3.1]{li2}. 
In fact, part of the proof of Theorem \ref{mainthm} (an argument in Section \ref{vK_section}) 
is inspired by \cite[\S 6.1]{li2}. 
However, we should emphasize that, we consider divisorial valuations over the 
\emph{original} $X$ in 
Theorem \ref{mainthm}, whereas \cite{li2} consider valuations over an 
\emph{affine cone} $(C_{r_0}, 0)$ of $X$, in particular $\dim C_{r_0}=\dim X+1$. 

\textbf{A postscript note:} 
While completing the current manuscript, the author learned that Chi Li 
independently showed 
essentially the same statement (\cite[Theorem 3.6]{li2}) as 
Corollary \ref{maincor} \eqref{maincor2} in the 
\emph{second} version of his preprint \cite{li2}. 
\end{enumerate}
\end{remark}

We explain the proof of Theorem \ref{mainthm}. It is known from \cite{B} and 
\cite{BBJ} (see Section \ref{stc_section}) that, uniform K-stability and K-semistability 
of a $\Q$-Fano variety $X$ is equivalent to \emph{uniform Ding stability} and 
\emph{Ding semistability} of $X$, respectively. In order to deduce the positivity of the 
function $\beta$ from uniform Ding stability or Ding semistability, we use an argument 
in \cite{fjt2}. Pick $r_0\in\Z_{>0}$ with $-r_0K_X$ Cartier. For any prime divisor 
$F$ over $X$, we consider the filtration $\sF$ of the graded $\Bbbk$-algebra 
$\bigoplus_{k\geq 0}H^0(X, \sO_X(-kr_0K_X))$ defined by 
\[
\sF^xH^0(X, \sO_X(-kr_0K_X)):=H^0(X, -kr_0K_X-xF)
\]
for $k\in\Z_{\geq 0}$ and $x\in\R_{\geq 0}$. We will construct a \emph{sequence} of 
test configurations for $(X, -K_X)$ from the filtration and calculate 
the corresponding Ding invariants. 
After taking a limit, we obtain the positivity of $\beta(F)$. 
For the converse, we heavily depend on the trick in \cite{LX} and the proof is 
inspired by \cite[\S 6.1]{li2}. It is enough to consider \emph{special test configurations} 
for $(X, -K_X)$ in order to test K-semistability and uniform K-stability 
(see \cite{LX} and Section \ref{stc_section}). Take any special test configuration 
$(\sX, -K_{\sX/\A^1})/\A^1$ for 
$(X, -K_X)$. We can show (see Section \ref{vK_section}) 
that the Donaldson-Futaki invariant $\DF(\sX, -K_{\sX/\A^1})$ of the 
special test configuration is a positive multiple of $\beta(v_{\sX_0})$, where 
$v_{\sX_0}$ is the \emph{restricted valuation} of the divisorial valuation on 
$X\times\A^1$ obtained by the fiber $\sX_0$ (see \cite{BHJ} and Section 
\ref{val_section}). Moreover the valuation $v_{\sX_0}$ is dreamy over $X$ 
(see \cite{BHJ, nystrom} and Section \ref{vK_section}). 
The positivity of $\DF(\sX, -K_{\sX/\A^1})$ immediately follows from 
the positivity of the function $\beta$.

This article is organized as follows. 
In Section \ref{tc_section}, we recall the notions of test configuration, 
K-stability, Ding stability and variants of these. 
Moreover, we recall their basic properties. 
In Section \ref{val_section}, we recall the theory of divisorial valuations and its 
restrictions. 
In Section \ref{filt_section}, we recall the theory of filtration of the graded linear series. 
In Section \ref{stc_section}, we show that it is enough to calculate the Donaldson-Futaki 
invariants of \emph{special} test configurations in order to test \emph{uniform} 
K-stability of a given $\Q$-Fano variety. We can prove the result 
in the same manner as \cite{LX, BBJ}. In Section \ref{Kv_section}, we show that 
uniform Ding stability (or Ding semistability) of a $\Q$-Fano variety implies 
the positivity of the function $\beta$. The proof is similar to the one in \cite{fjt2}. 
In Section \ref{vK_section}, we show the converse. Theorems \ref{mainthm} and 
\ref{mainK_thm} immediately 
follow from the results in Sections  \ref{stc_section}, 
\ref{Kv_section} and \ref{vK_section}. 
In Section \ref{log_section}, we consider a logarithmic version of Theorem \ref{mainthm}. 

\begin{ack}
The author thanks Doctor Chi Li, who sent him the preliminary version of 
\cite{li2}, Doctor Yuji Odaka, who gave him many comments, and Professors 
Robert Berman, Mattias Jonsson and Per Salberger, who gave him various 
suggestions during he visited Chalmers University of Technology. 
The author thanks the anonymous referees for their valuable comments 
and for a suggestion to add Section \ref{log_section}. 
The author was partially supported by a JSPS Fellowship for Young Scientists. 
\end{ack}

In this article, a \emph{variety} means an irreducible, reduced, 
separated scheme of finite type over $\Bbbk$. For the minimal model 
program, we refer the readers to \cite{KoMo}. 
We do not distinguish $\Q$-line bundles and $\Q$-Cartier divisors 
on varieties if there is no confusion. 
For a projective surjective morphism 
$\alpha\colon\sX\to C$ between normal varieties with $K_C$ $\Q$-Cartier, 
the \emph{relative canonical divisor} 
$K_{\sX/C}$ is defined by $K_{\sX/C}:=K_{\sX}-\alpha^*K_C$. For any point 
$t\in C$, the scheme-theoretic fiber $\alpha^{-1}(t)$ is denoted by $\sX_t$. 
For varieties $X_1$, $X_2$, the first (resp.\ the second) projection morphism 
$X_1\times X_2\to X_1$ (resp.\ $X_1\times X_2\to X_2$) is denoted by 
$p_1$ (resp.\ $p_2$).

\section{Preliminaries}\label{prelim_section}

\subsection{Test configurations}\label{tc_section}

We define the notions of test configuration, the Donaldson-Futaki invariant, 
the Ding invariant, and so on. For the notation, we basically follow 
\cite{BHJ}. For background, see \cite{tian1, don, RT} and references therein.

\begin{definition}\label{tc_dfn}
Let $X$ be a normal projective variety and $L$ be an ample $\Q$-line bundle. 
A \emph{test configuration} $(\sX, \sL)/\A^1$ for $(X, L)$ 
consists of the following data: 
\begin{itemize}
\item
a variety $\sX$ with a projective surjection $\alpha\colon\sX\to\A^1$, 
\item
a $\Q$-line bundle $\sL$ on $\sX$ which is nef over $\A^1$, 
\item
a $\G_m$-action on the pair $(\sX, \sL)$ such that the morphism $\alpha$ is 
$\G_m$-equivariant and $(\sX\setminus\sX_0, \sL|_{\sX\setminus\sX_0})$ is 
$\G_m$-equivariantly isomorphic to $(X\times(\A^1\setminus\{0\}), 
p_1^*L)$, where the action 
$\G_m\curvearrowright\A^1$ is defined multiplicatively. 
\end{itemize}
If $\sX$ (resp.\ $\sL$) is normal (resp.\ semiample, ample over $\A^1$), then we say 
that $(\sX, \sL)/\A^1$ is a \emph{normal} (resp.\ \emph{semiample, ample}) 
test configuration for $(X, L)$. 
\end{definition}

\begin{definition}\label{triv_dfn}
Let $X$ be a normal projective variety, let $L$ be an ample $\Q$-line bundle 
and let $(\sX, \sL)/\A^1$ be a test configuration for $(X, L)$. 
\begin{enumerate}
\renewcommand{\theenumi}{\arabic{enumi}}
\renewcommand{\labelenumi}{(\theenumi)}
\item\label{triv_dfn1}
We can glue $(\sX, \sL)/\A^1$ and $(X\times(\pr^1\setminus\{0\}), 
p_1^*L)$ along 
\[
(\sX\setminus\sX_0, \sL|_{\sX\setminus\sX_0})
\] 
and 
\[
(X\times
(\pr^1\setminus\{0,\infty\}), p_1^*L). 
\]
The compactification is denoted by $(\bar{\sX}, \bar{\sL})/\pr^1$. 
\item\label{triv_dfn2}
Assume that $(\sX, \sL)/\A^1$ is a semiample test configuration. Let $\nu\colon
\sX^\nu\to\sX$ be the normalization morphism and let $(\sX^\nu, \nu^*\sL)\to
(\sX^{\CAN}, \sL^{\CAN})$ be the canonical model of $\nu^*\sL$ over $\A^1$. 
(We note that $(\sX^{\CAN}, \sL^{\CAN})/\A^1$ is a normal, ample test configuration 
for $(X, L)$.) 
$(\sX, \sL)/\A^1$ is said to be a \emph{trivial} (resp.\ a \emph{product-type}) 
test configuration if $(\sX^{\CAN}, \sL^{\CAN})/\A^1$ is $\G_m$-equivariantly 
isomorphic to the pair $(X\times\A^1, p_1^*L+c\sX_0)$ for some $c\in\Q$ 
(resp.\  if $\sX^{\CAN}$ is isomorphic to $X\times\A^1$). 
\item\label{triv_dfn3}
Assume that $X$ is a $\Q$-Fano variety and $L=-K_X$. 
$(\sX, \sL)/\A^1$ is said to be a \emph{special test configuration} for $(X, -K_X)$ 
if it is a 
normal, ample test configuration for $(X, -K_X)$, $\sL=-K_{\sX/\A^1}$ and the pair 
$(\sX, \sX_0)$ is plt. 
\end{enumerate}
\end{definition}

\begin{definition}\label{DF_dfn}
Let $X$ be a normal projective variety of dimension $n$, $L$ be an ample 
$\Q$-line bundle and $(\sX, \sL)/\A^1$ be a 
test configuration for $(X, L)$. 
\begin{enumerate}
\renewcommand{\theenumi}{\arabic{enumi}}
\renewcommand{\labelenumi}{(\theenumi)}
\item\label{DF_dfn1}
Let 
\[\xymatrix{
& \bar{\sZ}  \ar[dl]_\Pi \ar[dr]^\Theta & \\
X\times\pr^1 & & \bar{\sX}
}\]
be the normalization of the graph of the birational map $X\times\pr^1\dashrightarrow
\bar{\sX}$. We set 
\[
\lambda_{\max}(\sX, \sL):=\frac{(\Pi^*p_1^*L^{\cdot n}\cdot\Theta^*\bar{\sL})}
{(L^{\cdot n})}.
\]
\item\label{DF_dfn2}
\cite{BHJ} We set 
\[
J^{\NA}(\sX, \sL):=\lambda_{\max}(\sX, \sL)-\frac{(\bar{\sL}^{\cdot n+1})}
{(n+1)(L^{\cdot n})}.
\]
\item\label{DF_dfn3}
\cite{B} (see also \cite[\S 7.7]{BHJ}) 
Assume that $X$ is a $\Q$-Fano variety and 
$(\sX, \sL)/\A^1$ is a normal test configuration for $(X, -K_X)$. 
Let $D_{(\sX, \sL)}$ be the $\Q$-divisor on $\bar{\sX}$ defined by 
\begin{itemize}
\item
$\Supp D_{(\sX, \sL)}\subset\sX_0$, 
\item
$D_{(\sX, \sL)}\sim_\Q -K_{\bar{\sX}/\pr^1}-\bar{\sL}$.
\end{itemize}
The \emph{Ding invariant} $\Ding(\sX, \sL)$ of $(\sX, \sL)/\A^1$ is defined by the 
following: 
\[
\Ding(\sX, \sL):=-\frac{(\bar{\sL}^{\cdot n+1})}{(n+1)((-K_X)^{\cdot n})}-1
+\lct(\sX, D_{(\sX, \sL)}; \sX_0),
\]
where
\[
\lct(\sX, D_{(\sX, \sL)}; \sX_0):=\max\{c\in\R\,|\,(\sX, D_{(\sX, \sL)}+c\sX_0)\text{: 
 sub log canonical}\}.
\]
\item\label{DF_dfn4}
\cite{wang,odk,LX}
Assume that $X$ is a $\Q$-Fano variety and 
$(\sX, \sL)/\A^1$ is a normal test configuration for $(X, -K_X)$. 
The \emph{Donaldson-Futaki invariant} $\DF(\sX, \sL)$ of $(\sX, \sL)/\A^1$ is 
defined by the following: 
\[
\DF(\sX, \sL):=\frac{n}{n+1}\cdot\frac{(\bar{\sL}^{\cdot n+1})}{((-K_X)^{\cdot n})}
+\frac{(\bar{\sL}^{\cdot n}\cdot K_{\bar{\sX}/\pr^1})}{((-K_X)^{\cdot n})}.
\]
\end{enumerate}
\end{definition}

\begin{remark}\label{rui_rmk}
Let $X$ be a normal projective variety of dimension $n$, $L$ be an ample $\Q$-line 
bundle and $(\sX, \sL)/\A^1$ be a 
semiample test configuration for $(X, L)$. The \emph{minimum norm} 
$\|(\sX, \sL)\|_m$ introduced by Dervan \cite{dervan} satisfies that 
\[
\frac{1}{n}\cdot J^{\NA}(\sX, \sL)\leq \frac{\|(\sX, \sL)\|_m}{(L^{\cdot n})}
\leq n\cdot J^{\NA}(\sX, \sL).
\]
See \cite[Remark 7.11 and Proposition 7.8]{BHJ}.
\end{remark}

\begin{proposition}\label{pull_prop}
Let $X$ be a $\Q$-Fano variety and $(\sX, \sL)/\A^1$ be a test 
configuration for $(X, -K_X)$. 
\begin{enumerate}
\renewcommand{\theenumi}{\arabic{enumi}}
\renewcommand{\labelenumi}{(\theenumi)}
\item\label{pull_prop1}
$($\cite[Theorem 1.3]{dervan} and \cite[Theorem 7.9]{BHJ}$)$
Assume that $(\sX, \sL)/\A^1$ is a normal, semiample test configuration. Then 
$J^{\NA}(\sX, \sL)\geq 0$ holds. Moreover, equality holds if and only if 
$(\sX, \sL)/\A^1$ is a trivial test configuration. 
\item\label{pull_prop2}
Let $\gamma\colon(\sY,\gamma^*\sL)\to(\sX, \sL)$ be a $\G_m$-equivariant 
projective birational morphism between test configurations for $(X, -K_X)$. 
Then the equality
$J^{\NA}(\sX, \sL)=J^{\NA}(\sY, \gamma^*\sL)$ holds. Moreover, if both $\sX$ and 
$\sY$ are normal, then $\Ding(\sX, \sL)=\Ding(\sY, \gamma^*\sL)$
and $\DF(\sX, \sL)=\DF(\sY, \gamma^*\sL)$ hold. 
\item\label{pull_prop3}
Assume that $(\sX, \sL)/\A^1$ is a normal test configuration. For any $d\in\Z_{>0}$, 
let $\psi_d\colon\sX^{(d)}\to \sX$ be the normalization of the base change of $\sX$ 
by the morphism $\A^1\to\A^1$ with $t\mapsto t^d$. Then the 
test configuration $(\sX^{(d)}, \psi_d^*\sL)/\A^1$ for $(X, -K_X)$ satisfies that 
$J^{\NA}(\sX^{(d)}, \psi_d^*\sL)=d\cdot J^{\NA}(\sX, \sL)$ and 
$\Ding(\sX^{(d)}, \psi_d^*\sL)=d\cdot \Ding(\sX, \sL)$. Moreover, if $\sL$ is ample 
over $\A^1$, then we have 
$\DF(\sX^{(d)}, \psi_d^*\sL)
\leq d\cdot\DF(\sX, \sL)$ and equality holds if and only if 
$\sX_0$ is reduced.
\item\label{pull_prop4}
$($\cite{B}$)$
Assume that $(\sX, \sL)/\A^1$ is a normal, ample test configuration. 
Then we have $\DF(\sX, \sL)\geq \Ding(\sX, \sL)$ and equality holds if and only if 
$\sL\sim_\Q -K_{\sX/\A^1}$ and the pair $(\sX, \sX_0)$ is log canonical. 
\end{enumerate}
\end{proposition}

\begin{proof}
For the proof of \eqref{pull_prop1}, \eqref{pull_prop2} and \eqref{pull_prop4}, 
see \cite[Theorem 1.3]{dervan}, \cite[Theorem 7.9]{BHJ}, \cite[\S 3]{B} and 
\cite[\S 3]{fjt2}. Let us prove \eqref{pull_prop3}. By the ramification formula 
(see \cite[p.\ 210]{LX}), we have 
$K_{\sX^{(d)}}+\sX^{(d)}_0=\psi_d^*(K_{\sX}+\red(\sX_0))$ and 
$K_{\bar{\sX}^{(d)}/\pr^1}=\psi_d^*(K_{\bar{\sX}/\pr^1}+\red(\sX_0)-\sX_0)$. 
Thus $D_{(\sX^{(d)}, \psi_d^*\sL)}=\psi^*(D_{(\sX, \sL)}-\red(\sX_0)+\sX_0)$. Thus we 
get 
$K_{\sX^{(d)}}+D_{(\sX^{(d)}, \psi_d^*\sL)}+(1-d+dc)\sX_0^{(d)}
=\psi_d^*(K_{\sX}+D_{(\sX, \sL)}
+c\sX_0)$ for any $c\in\Q$. Hence we have 
\[
-1+\lct(\sX^{(d)}, D_{(\sX^{(d)}, \psi_d^*\sL)}; \sX^{(d)}_0)=
d\cdot(-1+\lct(\sX, D_{(\sX, \sL)};\sX_0))
\]
(see \cite[Proposition 5.20]{KoMo}). This implies that 
$\Ding(\sX^{(d)}, \psi_d^*\sL)=d\cdot \Ding(\sX, \sL)$. 
The remaining assertions follow from \cite[Proposition 7.8]{BHJ} and \cite[Claim 1]{LX}. 
\end{proof}

\begin{remark}[{\cite[\S 7.3]{BHJ}}]\label{mabuchi_rmk}
Under the condition of Proposition \ref{pull_prop} \eqref{pull_prop3}, we have 
\[
M^{\NA}(\sX^{(d)}, \psi_d^*\sL)=d\cdot M^{\NA}(\sX, \sL), 
\]
where 
\[
M^{\NA}(\sX, \sL):=\DF(\sX, \sL)+\frac{1}{((-K_X)^{\cdot n})}
\left(\bar{\sL}^{\cdot n}\cdot \left(\red(\sX_0)-\sX_0\right)\right).
\]
\end{remark}

\begin{definition}\label{K_dfn}
Let $X$ be a $\Q$-Fano variety. 
\begin{enumerate}
\renewcommand{\theenumi}{\arabic{enumi}}
\renewcommand{\labelenumi}{(\theenumi)}
\item\label{K_dfn1}
\begin{enumerate}
\renewcommand{\theenumii}{\roman{enumii}}
\renewcommand{\labelenumii}{(\theenumii)}
\item\label{K_dfn11}
$X$ is said to be \emph{K-stable} (resp.\ \emph{K-semistable}) 
if for any nontrivial, normal, ample test configuration $(\sX, \sL)/\A^1$ for $(X, -K_X)$, 
we have $\DF(\sX, \sL)>0$ (resp.\ $\DF(\sX, \sL)\geq 0$).
\item\label{K_dfn12}
$X$ is said to be \emph{K-polystable} if 
$X$ is K-semistable and, if $\DF(\sX, \sL)=0$ 
for a normal, ample test configuration $(\sX, \sL)/\A^1$ for $(X, -K_X)$ implies that the 
configuration is a product-type. 
\item\label{K_dfn13}
$X$ is said to be \emph{uniformly K-stable} 
if there exists $\delta\in(0,1)$ such that 
$\DF(\sX, \sL)\geq \delta\cdot J^{\NA}(\sX, \sL)$ holds 
for any normal, ample test configuration $(\sX, \sL)/\A^1$ for $(X, -K_X)$. 
\end{enumerate}
\item\label{K_dfn2}
\begin{enumerate}
\renewcommand{\theenumii}{\roman{enumii}}
\renewcommand{\labelenumii}{(\theenumii)}
\item\label{K_dfn21}
$X$ is said to be \emph{Ding stable} (resp.\ \emph{Ding semistable}) 
if for any nontrivial, normal and 
ample test configuration $(\sX, \sL)/\A^1$ for $(X, -K_X)$, 
we have $\Ding(\sX, \sL)> 0$ (resp.\ $\Ding(\sX, \sL)\geq 0$). 
\item\label{K_dfn22}
$X$ is said to be \emph{Ding polystable} if 
$X$ is Ding semistable, and 
if $\Ding(\sX, \sL)=0$ for a normal, ample test configuration 
$(\sX, \sL)/\A^1$ for $(X, -K_X)$ with $\sX_0$ reduced implies that the configuration 
is a product-type. 
\item\label{K_dfn23}
$X$ is said to be \emph{uniformly Ding stable} if 
there exists $\delta\in(0,1)$ such that 
$\Ding(\sX, \sL)\geq \delta\cdot J^{\NA}(\sX, \sL)$ holds 
for any normal, ample test configuration $(\sX, \sL)/\A^1$ for $(X, -K_X)$. 
\end{enumerate}
\end{enumerate}
\end{definition}

\begin{remark}\label{KD_rmk}
\begin{enumerate}
\renewcommand{\theenumi}{\arabic{enumi}}
\renewcommand{\labelenumi}{(\theenumi)}
\item\label{KD_rmk1}
From Proposition \ref{pull_prop}, Ding semistability (resp.\ Ding polystability, 
Ding stability, uniform Ding stability) is stronger than K-semistability (resp.\ 
K-polystability, K-stability, uniform K-stability). It has been shown in \cite{BBJ} that 
Ding semistability (resp.\ uniform Ding stability) is equivalent to K-semistability
(resp.\ uniform K-stability). In fact, we will see in Section \ref{stc_section} that 
K-polystability (resp.\ K-stability) is also equivalent to Ding polystability (resp.\ 
Ding stability). 
\item\label{KD_rmk2}
The definition of Ding polystability in this article differs from the one in \cite{fjt2}. 
However, we will see in Section \ref{stc_section} that the definitions are equivalent. 
\end{enumerate}
\end{remark}

\subsection{Divisorial valuations}\label{val_section}

Let $X$ be a normal variety and let $K$ be the function field of $X$. 
We recall the results in \cite[\S 1, \S 4]{BHJ}. See also \cite{JM} and 
references therein. 
A \emph{divisorial valuation} on $X$ is a group homomorphism 
$v\colon K^*\to(\Q, +)$ of the form $c\cdot\ord_F$ with $c\in\Q_{>0}$ and $F$ 
a prime divisor over $X$.

\begin{definition}\label{beta_dfn}
Let $X$ be a $\Q$-Fano variety and 
let $v=c\cdot\ord_F$ be a divisorial valuation on $X$. 
For any $k\in\Z_{>0}$ with $-kK_X$ Cartier and for any $x\in\R_{\geq 0}$, we set
\[
H^0(X, -kK_X-xv):=H^0(X, -kK_X-cxF)\subset H^0(X, -kK_X).
\]
Similarly, we define $A_X(v):=c\cdot A_X(F)$, $\tau(v):=c^{-1}\cdot\tau(F)$, 
$\beta(v):=c\cdot\beta(F)$ and $j(v):=c\cdot j(F)$. 
$v$ is said to be \emph{dreamy} if $F$ is so. 
\end{definition}

From now on, let $X$ be a $\Q$-Fano variety with function field $K$, 
and $(\sX, \sL)/\A^1$ a normal test configuration for 
$(X, -K_X)$ such that there exists a 
projective birational $\G_m$-equivariant morphism $\Pi\colon\sX\to X\times\A^1$. 
Let $\sX_0=\sum_{i=0}^pm_iE_i$ be the irreducible decomposition of $\sX_0$ such that 
$E_0$ is the strict transform of $X\times\{0\}$ (thus $m_0=1$). 
Each $E_i$ gives a divisorial valuation 
\[
\ord_{E_i}\colon K(t)^*\to (\Q, +)
\]
on $X\times\A^1_t$. Let $r(\ord_{E_i})$ be the restriction of $\ord_{E_i}$ to $K^*$, 
let $v_{E_i}:=m_i^{-1}\cdot r(\ord_{E_i})$. Then $v_{E_0}\equiv 0$ and 
$v_{E_i}$ is a divisorial valuation on $X$ for any $1\leq i\leq p$ 
by \cite[Lemma 4.5]{BHJ}. 
We will use the following proposition later. 

\begin{proposition}[{\cite[Proposition 4.11]{BHJ}}]\label{discrep_prop}
Under the above notation, we have the equality 
\[
A_X(v_{E_i})=m_i^{-1}\cdot(\ord_{E_i}(K_{\sX/X\times\A^1})+1)-1
\]
for any $1\leq i\leq p$. 
\end{proposition}

\subsection{On filtrations}\label{filt_section}

We recall the theory of filtrations on graded linear series \cite{nystrom, BC, sz2, BHJ}. 

\begin{definition}[{see \cite[\S 1]{BC}}]\label{filt_dfn}
Let $X$ be a projective variety of dimension $n$, $L$ be an ample line bundle on $X$, 
$V_\bullet:=\bigoplus_{k\in\Z_{\geq 0}}V_k$ 
be the complete graded linear series 
of $L$, that is, $V_k=H^0(X, L^{\otimes k})$ for any $k\in\Z_{\geq 0}$. 
Let $\sF$ be a decreasing, left-continuous $\R$-filtration of the $\Bbbk$-algebra 
$V_\bullet$. 
\begin{enumerate}
\renewcommand{\theenumi}{\arabic{enumi}}
\renewcommand{\labelenumi}{(\theenumi)}
\item\label{filt_dfn1}
We say that $\sF$ is \emph{multiplicative} if 
$\sF^xV_k\otimes_\Bbbk\sF^{x'}V_{k'}$ maps to $\sF^{x+x'}V_{k+k'}$ for any 
$k$, $k'\in\Z_{\geq 0}$ and $x$, $x'\in\R$. 
\item\label{filt_dfn2}
We set 
\begin{eqnarray*}
e_{\max}(V_\bullet, \sF) & := & \limsup_{k\to\infty}\left(\frac{\sup\{x\in\R\,|\, 
\sF^x V_k\neq 0\}}{k}\right), \\
e_{\min}(V_\bullet, \sF) & := & \liminf_{k\to\infty}\left(\frac{\inf\{x\in\R\,|\, 
\sF^x V_k\neq V_k\}}{k}\right).
\end{eqnarray*}
We say that $\sF$ is \emph{linearly bounded} if both $e_{\max}(V_\bullet, \sF)$ and 
$e_{\min}(V_\bullet, \sF)$ are in $\R$. 
\item\label{filt_dfn3}
For any multiplicative $\sF$ and for any $x\in\R$, we define 
\[
\vol(\sF V_\bullet^x):=\limsup_{k\to\infty}\frac{\dim\sF^{kx}V_k}{k^n/n!}.
\]
\end{enumerate}
\end{definition}

We recall that test configurations induce filtrations. 

\begin{proposition}[{see \cite{nystrom, sz2, BHJ, odk15}}]\label{tcfilt_prop}
Let $X$ be a normal projective variety of dimension $n$, $L$ be an ample 
$\Q$-line bundle, $(\sX, \sL)/\A^1$ be a semiample 
test configuration for $(X, L)$ such that there exists a projective birational 
$\G_m$-equivariant morphism $\Pi\colon\sX\to X\times \A^1_t$, let 
$r_0\in\Z_{>0}$ with $r_0\sL$ Cartier, and $V_\bullet:=\bigoplus_{k\in\Z_{\geq 0}}V_k
=\bigoplus_{k\in\Z_{\geq 0}}H^0(X, \sO_X(kr_0L))$ be the complete graded linear 
series of $r_0L$. 
\begin{enumerate}
\renewcommand{\theenumi}{\arabic{enumi}}
\renewcommand{\labelenumi}{(\theenumi)}
\item\label{tcfilt_prop1}
We can define a decreasing, left-continuous $\R$-filtration $\sF_{(\sX, r_0\sL)}$ 
of $V_\bullet$ by 
\[
\sF^x_{(\sX, r_0\sL)}V_k:=\{f\in V_k\,|\, t^{\lfloor -x\rfloor}\Pi^*p_1^*f\in H^0(\sX, 
kr_0\sL)\}
\]
for $k\in\Z_{\geq 0}$ and $x\in\R$. 
Moreover, this filtration is multiplicative and linearly bounded. Furthermore, 
the graded $\Bbbk$-algebra 
\[
\bigoplus_{k\in\Z_{\geq 0},\,\,j\in\Z}\sF^j_{(\sX, r_0\sL)} V_k
\]
is finitely generated. 
\item\label{tcfilt_prop2}
For any $k\in\Z_{\geq 0}$, we set 
\begin{eqnarray*}
\lambda_{\max}^{(k)} & := & \sup\{x\in\R\,|\, \sF^x_{(\sX, r_0\sL)}V_k\neq 0\}, \\
\lambda_{\min}^{(k)} & := & \inf\{x\in\R\,|\, \sF^x_{(\sX, r_0\sL)}V_k\neq V_k\}, \\
w(k) & := & \int_{\lambda_{\min}^{(k)}}^{\lambda_{\max}^{(k)}}\dim\sF^x_{(\sX, r_0\sL)}V_k
dx+\lambda_{\min}^{(k)}\cdot\dim V_k.
\end{eqnarray*}
\begin{enumerate}
\renewcommand{\theenumii}{\roman{enumii}}
\renewcommand{\labelenumii}{(\theenumii)}
\item\label{tcfilt_prop21}
We have
\[
\lambda_{\max}(\sX, \sL)=\sup_{k\to\infty}\frac{\lambda_{\max}^{(k)}}{kr_0}
=\lim_{k\to\infty}\frac{\lambda_{\max}^{(k)}}{kr_0}.
\]
\item\label{tcfilt_prop22}
$w(k)$ is a polynomial of degree at most $n+1$ for $k\gg 0$. Moreover, 
we have 
\[
\lim_{k\to\infty}\frac{w(k)}{k\cdot \dim V_k}=\frac{r_0\cdot(\bar{\sL}^{\cdot n+1})}
{(n+1)(L^{\cdot n})}.
\]
In particular, we have 
\[
\lim_{k\to\infty}\frac{w(k)}{k^{n+1}/(n+1)!}=r_0^{n+1}(\bar{\sL}^{\cdot n+1}).
\]
\item\label{tcfilt_prop23}
Assume that $X$ is a $\Q$-Fano variety and $L=-K_X$. 
Let us consider the asymptotic expansion
\[
\frac{w(k)}{k\cdot\dim V_k}=F_0+k^{-1}F_1+k^{-2}F_2+\cdots.
\]
If $\sX$ is normal, then we have $\DF(\sX, \sL)=-2 F_1$. 
\item\label{tcfilt_prop24}
\cite{odk15} Assume that $X$ is a $\Q$-Fano variety and $L=-K_X$. 
Let us set $\DF(\sX, \sL):=-2F_1$ for possibly non-normal $\sX$. 
If $(X, -K_X)$ is K-stable $($resp.\ K-semistable$)$, $(\sX, \sL)/\A^1$ is nontrivial, 
and $\Pi$ is not an isomorphism in 
codimension one, then $\DF(\sX, \sL)>0$ $($resp.\ $\geq 0)$ holds. 
\end{enumerate}
\end{enumerate}
\end{proposition}

\begin{proof}
For \eqref{tcfilt_prop1}, see \cite[\S 2.5]{BHJ} and \cite{nystrom} for example. 
For \eqref{tcfilt_prop2}, see \cite[\S 5, Lemma 7.7, Proposition 3.12]{BHJ} and 
\cite{odk15} for example. 
\end{proof}

\section{K-stability and Ding stability}\label{stc_section}

In this section, we recall the results \cite{LX, BBJ}. Moreover, we see that 
it is enough to consider special test configurations in order to test uniform K-stability 
of $\Q$-Fano varieties. 
Many results in this section are already known. The author wrote down this section 
in detail just for the readers' convenience. 

\begin{thm}[{see \cite[Theorem 2]{LX}}]\label{LX2_thm}
Let $X$ be a $\Q$-Fano variety and $(\sX, \sL)/\A^1$ be a normal, ample 
test configuration for $(X, -K_X)$. Then there exist $d\in\Z_{>0}$, a projective birational 
$\G_m$-equivariant morphism $\pi\colon \sX^{\LC}\to\sX^{(d)}$ $($with $\psi_d
\colon\sX^{(d)}\to\sX$ as in Proposition \ref{pull_prop} \eqref{pull_prop3}$)$, and 
a normal, ample test configuration $(\sX^{\LC}, \sL^{\LC})/\A^1$ for $(X, -K_X)$ 
such that: 
\begin{enumerate}
\renewcommand{\theenumi}{\arabic{enumi}}
\renewcommand{\labelenumi}{(\theenumi)}
\item\label{LX2_thm1}
$(\sX^{\LC}, \sX_0^{\LC})$ is log canonical. 
\item\label{LX2_thm2}
$\DF(\sX^{\LC}, \sL^{\LC})\leq d\cdot\DF(\sX, \sL)$ holds. Moreover, 
equality holds if and only if $(\sX, \sX_0)$ is log canonical and $-K_{\sX^{(d)}}\sim_\Q 
\psi_d^*\sL$; in which case $\sL^{\LC}\sim_\Q\pi^*\psi_d^*\sL$ and $\pi$
is an isomorphism. 
\item\label{LX2_thm3}
$\Ding(\sX^{\LC}, \sL^{\LC})\leq d\cdot\Ding(\sX, \sL)$ holds. Moreover, 
equality holds if and only if $(\sX^{(d)}, \sX_0^{(d)})$ is log canonical and 
$-K_{\sX^{(d)}}\sim_\Q\psi_d^*\sL$; 
in which case $\sL^{\LC}\sim_\Q\pi^*\psi_d^*\sL$ and $\pi$
is an isomorphism. 
\item\label{LX2_thm4}
For any $\delta\in[0,1]$, we have the inequality 
\[
d(\Ding(\sX, \sL)-\delta\cdot J^{\NA}(\sX, \sL))\geq
\Ding(\sX^{\LC}, \sL^{\LC})-\delta\cdot J^{\NA}(\sX^{\LC}, \sL^{\LC}).
\]
\end{enumerate}
\end{thm}

\begin{proof}
We repeat the proof of \cite[Theorem 2]{LX}. 
By the semistable reduction theorem and \cite[Proposition 2]{LX}, there exist 
$d\in\Z_{>0}$ and the log canonical modification $\pi\colon\sX^{\LC}\to(\sX^{(d)}, 
\sX_0^{(d)})$ of the pair $(\sX^{(d)}, \sX_0^{(d)})$. Set $\sL_0^{\LC}:=\pi^*\psi_d^*\sL$. 
Let $E$ be the $\Q$-divisor on $\bar{\sX}^{\LC}$ defined by 
\begin{itemize}
\item
$\Supp E\subset\sX_0^{\LC}$, 
\item
$E\sim_\Q K_{\bar{\sX}^{\LC}/\pr^1}+\bar{\sL}_0^{\LC}$.
\end{itemize}
Set $\sL_t^{\LC}:=\sL_0^{\LC}+tE$. Then, by \cite[Theorem 2 and the proof of 
Proposition 3]{LX}, 
$(\sX^{\LC}, \sL_t^{\LC})/\A^1$ is a normal, ample test configuration for 
$(X, -K_X)$ and 
satisfies the condition \eqref{LX2_thm2} for any $0< t\ll 1$. 

We check \eqref{LX2_thm3}. Let $\sX_0^{\LC}=\sum_{i=1}^pE_i$ be the irreducible 
decomposition and set $E=:\sum_{i=1}^pe_iE_i$. We can assume that $e_1\leq\dots\leq
e_p$. Since $-K_{\bar{\sX}^{\LC}/\pr^1}-\bar{\sL}_t^{\LC}\sim_\Q-(1+t)E$, we have 
$D_{(\sX^{\LC}, \sL_t^{\LC})}=-(1+t)E$. Thus 
\[
\lct(\sX^{\LC}, D_{(\sX^{\LC}, \sL_t^{\LC})}; \sX_0^{\LC})=1+(1+t)e_1
\]
since $(\sX^{\LC}, \sum_{i=1}^pE_i)$ is log canonical. Hence we have 
\[
\Ding(\sX^{\LC}, \sL_t^{\LC})=-\frac{(\bar{\sL}_t^{\LC}{}^{\cdot n+1})}
{(n+1)((-K_X)^{\cdot n})}+(1+t)e_1.
\]
This implies that 
\begin{eqnarray*}
&&\frac{d}{dt}\Ding(\sX^{\LC}, \sL_t^{\LC})=-\frac{(\bar{\sL}_t^{\LC}{}^{\cdot n}
\cdot E)}{((-K_X)^{\cdot n})}+e_1=
\frac{\left(\bar{\sL}_t^{\LC \cdot n}\cdot
\left(e_1\sX_0^{\LC}-E\right)\right)}{((-K_X)^{\cdot n})}\\
&=&\frac{1}{((-K_X)^{\cdot n})}\left(\bar{\sL}_t^{\LC \cdot n}\cdot
\sum_{i=1}^p(e_1-e_i)E_i\right)\leq 0
\end{eqnarray*}
for any $0< t\ll 1$. Thus we get 
\[
\Ding(\sX^{\LC}, \sL^{\LC}_t)\leq\Ding(\sX^{\LC}, \sL_0^{\LC})=d\cdot\Ding(\sX, \sL)
\]
(see Proposition \ref{pull_prop} \eqref{pull_prop3}). If the above inequality is an 
equality, then 
$E\sim_{\Q, \pr^1} 0$ since $e_1=\cdots=e_p$. This implies that $\pi$ is an isomorphism
since $E$ is $\pi$-ample. 

We check \eqref{LX2_thm4}. Let 
\[\xymatrix{
& \bar{\sZ}  \ar[dl]_\Pi \ar[dr]^\Theta & \\
X\times\pr^1 & & \bar{\sX}^{\LC}
}\]
be the normalization of the graph. Set $\phi_t:=\Theta^*\bar{\sL}_t^{\LC}$, 
$\phi_{\triv}:=\Pi^*p_1^*(-K_X)$, $V:=((-K_X)^{\cdot n})$ as in \cite[\S 6.1]{BHJ}. Then 
\begin{eqnarray*}
&&(n+1)V(d(\Ding(\sX, \sL)-\delta J^{\NA}(\sX, \sL))\\
&&-(\Ding(\sX^{\LC}, \sL_t^{\LC})-
\delta J^{\NA}(\sX^{\LC}, \sL_t^{\LC})))\\
&=&(n+1)V(\Ding(\sX^{\LC}, \sL_0^{\LC})-\delta J^{\NA}(\sX^{\LC}, \sL_0^{\LC})\\
&&-(\Ding(\sX^{\LC}, \sL_t^{\LC})-
\delta J^{\NA}(\sX^{\LC}, \sL_t^{\LC})))\\
&=&(1-\delta)((\phi_t^{\cdot n+1})-(\phi_0^{\cdot n+1}))
+\delta(n+1)t(\phi_{\triv}^{\cdot n}\cdot\Theta^*E)-(n+1)te_1V\\
&=&(1-\delta)t\sum_{i=0}^n\left((\phi_t^{\cdot i}\cdot
\phi_0^{\cdot n-i}\cdot\Theta^*E)-e_1V\right)\\
&&+\delta(n+1)t\left((\phi_{\triv}^{\cdot n}\cdot\Theta^*E)-e_1V\right)\\
&=&(1-\delta)t\sum_{i=0}^n\left(\phi_t^{\cdot i}\cdot\phi_0^{\cdot n-i}\cdot
\Theta^*\sum_{j=1}^p(e_j-e_1)E_j\right)\\
&&+\delta(n+1)t\left(\phi_{\triv}^{\cdot n}\cdot\Theta^*\sum_{j=1}^p
(e_j-e_1)E_j\right)
\geq 0.
\end{eqnarray*}
Thus we get the desired inequality. 
\end{proof}

\begin{thm}[{see \cite[Theorem 3]{LX} and \cite[Theorem 2.1]{BBJ}}]\label{LX3_thm}
Let $X$ be a $\Q$-Fano variety and $(\sX, \sL)/\A^1$ be a normal, ample 
test configuration for $(X, -K_X)$ with $(\sX, \sX_0)$ log canonical. 
Then there exists 
a normal, ample test configuration $(\sX^{\AC}, \sL^{\AC})/\A^1$ for $(X, -K_X)$ 
with $(\sX^{\AC}, \sX^{\AC}_0)$ log canonical such that: 
\begin{enumerate}
\renewcommand{\theenumi}{\arabic{enumi}}
\renewcommand{\labelenumi}{(\theenumi)}
\item\label{LX3_thm1}
$-K_{\sX^{\AC}}\sim_\Q\sL^{\AC}$. 
\item\label{LX3_thm2}
$\DF(\sX^{\AC}, \sL^{\AC})\leq \DF(\sX, \sL)$ holds. Moreover, 
equality holds if and only if $(\sX^{\AC}, \sL^{\AC})\simeq(\sX, \sL)$ over $\A^1$. 
\item\label{LX3_thm3}
$\Ding(\sX^{\AC}, \sL^{\AC})\leq \Ding(\sX, \sL)$ holds. Moreover, 
equality holds if and only if $(\sX^{\AC}, \sL^{\AC})\simeq(\sX, \sL)$ over $\A^1$. 
\item\label{LX3_thm4}
For any $\delta\in[0,1]$, we have the inequality 
\[
\Ding(\sX, \sL)-\delta\cdot J^{\NA}(\sX, \sL)\geq
\Ding(\sX^{\AC}, \sL^{\AC})-\delta\cdot J^{\NA}(\sX^{\AC}, \sL^{\AC}).
\]
\end{enumerate}
\end{thm}

\begin{proof}
We repeat the proof of \cite[Theorem 3]{LX}. 
Since $(\sX, \sX_0)$ is log canonical, $\sX$ is log terminal. Thus there exists 
a $\G_m$-equivariant small $\Q$-factorial modification $\sigma\colon\sX^0\to\sX$ 
by \cite[Corollary 1.4.3]{BCHM}. Set $\sL^0:=\sigma^*\sL$. Since $\sL$ is ample 
over $\A^1$, $-K_{\sX^0/\A^1}+l\sL^0$ is semiample and big over $\A^1$ for $l\gg 1$. 
Thus 
\[
\sH^0:=\sL^0-\frac{1}{l+1}(\sL^0+K_{\sX^0/\A^1})=\frac{1}{l+1}(l\sL^0-K_{\sX^0/\A^1})
\]
is nef over $\A^1$. As in \cite[p.\ 211]{LX}, we run the $K_{\sX^0/\A^1}$-MMP 
over $\A^1$ with scaling $\sH^0$. 
More precisely, we set $\lambda_0:=l+1$, 
\[
\lambda_{j+1}:=\min\{\lambda\,|\, K_{\sX^j/\A^1}+\lambda\sH^j
\text{ is nef over }\A^1\},
\]
$\sX^j\dashrightarrow\sX^{j+1}/\A^1$ the $\G_m$-equivariant birational map 
corresponds to a $K_{\sX^j/\A^1}$-negative and 
$(K_{\sX^j/\A^1}+\lambda_{j+1}\sH^j)$-trivial extremal ray, and let $\sH^{j+1}$ be the 
strict transform of $\sH^j$ on $\sX^{j+1}$. Then we get 
\[
\sX^0\dashrightarrow\sX^1\dashrightarrow\cdots\dashrightarrow\sX^k 
\]
and $l+1=\lambda_0>\lambda_1\geq\cdots\geq\lambda_k>\lambda_{k+1}=1$. 
For any $0\leq j\leq k-1$ and $\lambda\in[\lambda_{j+1}, \lambda_j]$, we set 
\[
\sL_\lambda^j:=\frac{1}{\lambda-1}(K_{\sX^j/\A^1}+\lambda\sH^j). 
\]
The pair $(\sX^j, \sL_\lambda^j)/\A^1$ is not a test configuration for $(X, -K_X)$ 
in general 
because of the small $\Q$-factorial modification $\sigma$. However, we can still define 
$\DF(\sX^j, \sL_\lambda^j)$, $\Ding(\sX^j, \sL_\lambda^j)$ and 
$J^{\NA}(\sX^j, \sL_\lambda^j)$ in the same way. (In fact, $\sL_\lambda^j$ is 
semiample over $\A^1$ by the base point free theorem and the canonical model of 
$\sL_\lambda^j$ over $\A^1$ gives a normal, ample test configuration for $(X, -K_X)$.)
By \cite[Lemma 2]{LX}, the canonical model $\mu\colon (\sX^k, \sL_{\lambda_k}^k)
\to(\sX^{\AC}, \sL^{\AC})$ of $\sL_{\lambda_k}^k$ over $\A^1$ satisfies 
$\sL^{\AC}\sim_\Q-K_{\sX^{\AC}}$. Moreover, by \cite[Theorem 3]{LX}, 
$(\sX^{\AC}, \sL^{\AC})/\A^1$ satisfies condition \eqref{LX3_thm2}. 
Let $E$ be the $\Q$-divisor on $\bar{\sX}^0$ defined by 
\begin{itemize}
\item
$\Supp E\subset\sX_0^0$, 
\item
$E\sim_\Q K_{\bar{\sX}^0/\pr^1}+\bar{\sH}^0$. 
\end{itemize}
Since $-(K_{\bar{\sX}^j/\pr^1}+\bar{\sL}_\lambda^j)=-(\lambda/(\lambda-1))
(K_{\bar{\sX}^j/\pr^1}+\bar{\sH}^j)$, we have 
$D_{(\sX^j, \sL_\lambda^j)}=-(\lambda/(\lambda-1))E^j$, where $E^j$ is the strict 
transform of $E$ on $\sX^j$. 
By Proposition \ref{pull_prop}, we have 
$\Ding(\sX, \sL)=\Ding(\sX^0, \sL^0)$, $\Ding(\sX^j, \sL^j_{\lambda_{j+1}})
=\Ding(\sX^{j+1}, \sL^{j+1}_{\lambda_{j+1}})$, $\Ding(\sX^k, \sL^k_{\lambda_k})=
\Ding(\sX^{\AC}, \sL^{\AC})$ and 
$J^{\NA}(\sX, \sL)=J^{\NA}(\sX^0, \sL^0)$, $J^{\NA}(\sX^j, \sL^j_{\lambda_{j+1}})=
J^{\NA}(\sX^{j+1}, \sL^{j+1}_{\lambda_{j+1}})$ and $J^{\NA}(\sX^k, \sL^k_{\lambda_k})=
J^{\NA}(\sX^{\AC}, \sL^{\AC})$ for any $0\leq j\leq k-1$. 

We check \eqref{LX3_thm3}. For any $0\leq j\leq k-1$, let $\sX_0^j=\sum_{i=1}^pE_i$ be 
the irreducible decomposition and set $E^j=\sum_{i=1}^pe_iE_i.$ We can assume that 
$e_1\leq\cdots\leq e_p$. Take any $\lambda\in[\lambda_{j+1}, \lambda_j]$. 
We note that 
\[
\Ding(\sX^j, \sL^j_\lambda)=-\frac{(\bar{\sL}_\lambda^{j}{}^{\cdot n+1})}
{(n+1)((-K_X)^{\cdot n})}+\frac{\lambda}{\lambda-1}e_1.
\]
Thus 
\begin{eqnarray*}
&&\frac{d}{d\lambda}\Ding(\sX^j, \sL^j_\lambda)=\frac{1}{(\lambda-1)^2}\left(
\frac{(\bar{\sL}_\lambda^j{}^{\cdot n}\cdot E^j)}{((-K_X)^{\cdot n})}-e_1\right)\\
&=&\frac{1}{(\lambda-1)^2((-K_X)^{\cdot n})}\left(\bar{\sL}_{\lambda}^{j\cdot n}
\cdot\sum_{i=1}^p(e_i-e_1)E_i\right)
\geq 0.
\end{eqnarray*}
This implies that $\Ding(\sX, \sL)\geq \Ding(\sX^{\AC}, \sL^{\AC})$. Assume that the 
inequality is an equality. Then $\Ding(\sX^0, \sL_{\lambda_0}^0)=\Ding(\sX^0, 
\sL_{\lambda_1}^0)$. We know that $\lambda_0>\lambda_1$ since $l\gg 0$. 
Moreover, since $\sigma$ is small, we have 
$(\bar{\sL}_\lambda^0{}^{\cdot n}\cdot E_i)>0$ for any $\lambda_1\ll\lambda\leq 
\lambda_0$ and any irreducible component $E_i\subset\sX_0^0$. Thus we get 
$E\sim_{\Q,\pr^1}0$. This implies that $-K_{\sX^0}\sim_\Q\sL^0$. 
In this case, $k=0$ and $(\sX^{\AC}, \sL^{\AC})\simeq(\sX, \sL)$ over $\A^1$ hold. 

We check \eqref{LX3_thm4}. Take any $0\leq j\leq k-1$. 
Let 
\[\xymatrix{
& \bar{\sZ}  \ar[dl]_\Pi \ar[dr]^\Theta & \\
X\times\pr^1 & & \bar{\sX}^j
}\]
be the normalization of the graph. Set $\phi_\lambda:=\Theta^*\bar{\sL}_\lambda^j$, 
$\phi_{\triv}:=\Pi^*p_1^*(-K_X)$, $V:=((-K_X)^{\cdot n})$ as in \cite[\S 6.1]{BHJ}. Then 
\[
\phi_{\lambda_{j+1}}-\phi_{\lambda_j}=
\left(\frac{\lambda_{j+1}}{\lambda_{j+1}-1}-\frac{\lambda_j}{\lambda_j-1}\right)
\Theta^*E^j.
\]
Thus we have
\begin{eqnarray*}
&&(n+1)V(\Ding(\sX^j, \sL_{\lambda_j}^j)-\delta J^{\NA}(\sX^j, \sL_{\lambda_j}^j)\\
&&-(\Ding(\sX^j, \sL_{\lambda_{j+1}}^j)-
\delta J^{\NA}(\sX^j, \sL_{\lambda_{j+1}}^j)))\\
&=&(1-\delta)\left(\frac{\lambda_{j+1}}{\lambda_{j+1}-1}-\frac{\lambda_j}
{\lambda_j-1}\right)\sum_{i=0}^n\left((\phi_{\lambda_{j+1}}^{\cdot i}\cdot
\phi_{\lambda_j}^{\cdot n-i}\cdot\Theta^*E^j)-e_1V\right) \\
&&+\delta(n+1)\left(\frac{\lambda_{j+1}}{\lambda_{j+1}-1}-\frac{\lambda_j}
{\lambda_j-1}\right)\left((\phi_{\triv}^{\cdot n}\cdot\Theta^*E^j)-e_1V\right)\geq 0.
\end{eqnarray*}
Therefore we get the desired inequality. 
\end{proof}

\begin{thm}[{see \cite[Theorem 4]{LX}}]\label{LX4_thm}
Let $X$ be a $\Q$-Fano variety and $(\sX, -K_{\sX/\A^1})/\A^1$ be a normal, ample 
test configuration for $(X, -K_X)$ with $(\sX, \sX_0)$ log canonical. 
Then there exist $d\in\Z_{>0}$ and a projective birational 
$\G_m$-equivariant birational map $\sX^{(d)}
\dashrightarrow\sX^{\s}$ over $\A^1$ $($with $\psi_d
\colon\sX^{(d)}\to\sX$ as in Proposition \ref{pull_prop} \eqref{pull_prop3}$)$ 
such that: 
\begin{enumerate}
\renewcommand{\theenumi}{\arabic{enumi}}
\renewcommand{\labelenumi}{(\theenumi)}
\item\label{LX4_thm1}
$(\sX^{\s}, -K_{\sX^{\s}/\A^1})/\A^1$ is a special test configuration for $(X, -K_X)$. 
\item\label{LX4_thm2}
$\DF(\sX^{\s}, -K_{\sX^{\s}/\A^1})\leq d\cdot\DF(\sX, -K_{\sX/\A^1})$ holds. 
Moreover, 
equality holds if and only if $(\sX, -K_{\sX/\A^1})/\A^1$ is a special test configuration 
for $(X, -K_X)$. We remark that, by Proposition 
\ref{pull_prop} \eqref{pull_prop4}, we have the equalities
$\DF(\sX^{\s}, -K_{\sX^{\s}/\A^1})=\Ding(\sX^{\s}, 
-K_{\sX^{\s}/\A^1})$ and $\DF(\sX, -K_{\sX/\A^1})=\Ding(\sX, -K_{\sX/\A^1})$. 
\item\label{LX4_thm3}
For any $\delta\in[0,1]$, we have the inequality 
\begin{eqnarray*}
&&d(\Ding(\sX, -K_{\sX/\A^1})-\delta\cdot J^{\NA}(\sX, -K_{\sX/\A^1}))\\
&\geq&\Ding(\sX^{\s}, -K_{\sX^{\s}/\A^1})
-\delta\cdot J^{\NA}(\sX^{\s}, -K_{\sX^{\s}/\A^1}).
\end{eqnarray*}
\end{enumerate}
\end{thm}

\begin{proof}
We repeat the proof of \cite[Theorem 4]{LX}. 
By \cite[Theorem 4]{LX}, there exist $d\in\Z_{>0}$ and a $\G_m$-equivariant birational 
map $\sX^{(d)}\dashrightarrow\sX^{\s}$ such that the conditions \eqref{LX4_thm1} and 
\eqref{LX4_thm2} hold and the discrepancy $a(\sX_0^{\s}; \sX^{(d)}, 0)$ is equal to zero. 
By \cite[Corollary 1.4.3]{BCHM}, there exists a $\G_m$-equivariant birational morphism 
$\pi'\colon\sX'\to\sX^{(d)}$ with $\sX'$ normal and $\Q$-factorial such that 
$\pi'$ exactly contracts the divisor $\sX_0^{\s}$. 
We know that $K_{\sX'}={\pi'}^*K_{\sX^{(d)}}$. 
Thus we have $\Ding(\sX', -K_{\sX'/\A^1})=\Ding(\sX^{(d)}, -K_{\sX^{(d)}/\A^1})=d\cdot
\Ding(\sX, -K_{\sX/\A^1})$ and $J^{\NA}(\sX', -K_{\sX'/\A^1})=J^{\NA}
(\sX^{(d)}, -K_{\sX^{(d)}/\A^1})=d\cdot J^{\NA}(\sX, -K_{\sX/\A^1})$. 
Consider a common partial resolution 
\[\xymatrix{
X\times\pr^1 & \bar{\sZ}  \ar[l]_\Pi \ar[d]^\Theta \ar[r]^\Xi & \bar{\sX}^{\s}.\\
 & \bar{\sX}' &
}\]
Set $E:=\Xi^*(-K_{\bar{\sX}^{\s}/\pr^1})-\Theta^*(-K_{\bar{\sX}'/\pr^1})$. 
Then $-E$ is $\Xi$-exceptional and $\Xi$-nef. Thus $E$ is effective by 
the negativity lemma. Set $\phi_0:=\Theta^*(-K_{\bar{\sX}'/\pr^1})$, $\phi_1:=
\Xi^*(-K_{\bar{\sX}^{\s}/\pr^1})$, $\phi_{\triv}:=\Pi^*p_1^*(-K_X)$ and 
$V:=((-K_X)^{\cdot n})$ as in \cite[\S 6.1]{BHJ}. Then we have 
\begin{eqnarray*}
&&(n+1)V(\Ding(\sX', -K_{\sX'/\A^1})-\delta J^{\NA}(\sX', -K_{\sX'/\A^1})\\
&&-(\Ding(\sX^{\s}, -K_{\sX^{\s}/\A^1})-
\delta J^{\NA}(\sX^{\s}, -K_{\sX^{\s}/\A^1})))\\
&=&(1-\delta)\sum_{i=0}^n(\phi_1^{\cdot i}\cdot\phi_0^{\cdot n-i}\cdot E)
+\delta(n+1)(\phi_{\triv}^{\cdot n}\cdot E)\geq 0.
\end{eqnarray*}
Therefore we get the desired inequality. 
\end{proof}

By Theorems \ref{LX2_thm}, \ref{LX3_thm} and \ref{LX4_thm}, we immediately get 
the following corollary: 

\begin{corollary}\label{stc_cor}
Let $X$ be a $\Q$-Fano variety. 
\begin{enumerate}
\renewcommand{\theenumi}{\arabic{enumi}}
\renewcommand{\labelenumi}{(\theenumi)}
\item\label{stc_cor1}
$($see \cite[Theorem 2.1]{BBJ} and \cite[Corollary 1]{LX}$)$
For any $\delta\in[0,1)$, the following conditions are equivalent: 
\begin{enumerate}
\renewcommand{\theenumii}{\roman{enumii}}
\renewcommand{\labelenumii}{(\theenumii)}
\item\label{stc_cor11}
$\DF(\sX, \sL)\geq \delta\cdot J^{\NA}(\sX, \sL)$ holds for any normal, ample 
test configuration $(\sX, \sL)/\A^1$ for $(X, -K_X)$. 
\item\label{stc_cor12}
$\Ding(\sX, \sL)\geq \delta\cdot J^{\NA}(\sX, \sL)$ holds for any normal, ample 
test configuration $(\sX, \sL)/\A^1$ for $(X, -K_X)$. 
\item\label{stc_cor13}
$\DF(\sX, -K_{\sX/\A^1})\geq \delta\cdot J^{\NA}(\sX, -K_{\sX/\A^1})$ holds 
for any special test configuration $(\sX, -K_{\sX/\A^1})/\A^1$ for $(X, -K_X)$. 
\end{enumerate}
\item\label{stc_cor2}
$X$ is K-stable $($resp.\ K-polystable$)$ if and only if 
$X$ is Ding stable $($resp.\ Ding polystable$)$. 
\end{enumerate}
\end{corollary}

\section{Sequences of test configurations}\label{Kv_section}

In this section, we prove the following theorem. 

\begin{thm}\label{Kv_thm}
Let $X$ be a $\Q$-Fano variety. Assume that there exists $\delta\in[0,1)$ such that 
$\Ding(\sX, \sL)\geq \delta \cdot J^{\NA}(\sX, \sL)$ holds for any 
normal, semiample test configuration $(\sX, \sL)/\A^1$ for $(X, -K_X)$. Then 
$\beta(F)\geq\delta\cdot j(F)$ holds for any prime divisor $F$ over $X$. 
\end{thm}

\begin{proof}
Take any projective log resolution $\sigma\colon Y\to X$ with $F$ a smooth divisor 
on $Y$. Fix $r_0\in\Z_{>0}$ with $-r_0K_X$ Cartier, and let $V_\bullet=
\bigoplus_{k\in\Z_{\geq 0}}V_k$ be the complete graded linear series of $-r_0K_X$. 
Let $\sF$ be the decreasing, left-continuous $\R$-filtration of $V_\bullet$ defined by 
\[
\sF^xV_k:=\begin{cases}
H^0(X, -kr_0K_X-xF) & \text{if } x\in\R_{\geq 0},\\
V_k & \text{otherwise,}
\end{cases}
\]
for $k\in\Z_{\geq 0}$ and $x\in\R$. 
Then $\sF$ is multiplicative and linearly bounded. In fact, we have 
$e_{\max}(V_\bullet,\sF)=r_0\cdot\tau(F)$ and $e_{\min}(V_\bullet, \sF)=0$. 
Set 
\[
I_{(k,x)}:=\Image(\sF^xV_k\otimes_\Bbbk\sO_X(kr_0K_X)\to\sO_X)
\]
for $k\in\Z_{\geq 0}$ and $x\in\R$, where the homomorphism is the evaluation. 
By definition, $I_{(k,x)}\cdot\sO_Y$ is equal to 
\[
\Image\left(H^0\left(Y, 
\sigma^*\sO_X(-kr_0K_X)(\lfloor -xF\rfloor)\right)\otimes_\Bbbk
\sigma^*\sO_X(kr_0K_X)\to\sO_Y\right).
\]
In particular, $I_{(k,x)}\cdot\sO_Y\subset\sO_Y(\lfloor -xF\rfloor)$ holds. 

\begin{claim}\label{sat_claim}
For any $k\in\Z_{\geq 0}$ and $x\in\R$, we have the equality
\[
\sF^xV_k=H^0(X, \sO_X(-kr_0K_X)\cdot I_{(k,x)}). 
\]
In other words, the filtration $\sF$ is 
\emph{saturated} in the sense of \cite[Definition 4.4]{fjt2}. 
\end{claim}

\begin{proof}[Proof of Claim \ref{sat_claim}]
By \cite[Proposition 4.3 (5)]{fjt2}, 
\[
\sF^xV_k\subset H^0(X, \sO_X(-kr_0K_X)\cdot I_{(k,x)})
\]
holds. On the other hand, $H^0(X, \sO_X(-kr_0K_X)\cdot I_{(k,x)})$ maps to 
\begin{eqnarray*}
&&H^0\left(Y, \sigma^*\left(\sO_X(-kr_0K_X)\otimes I_{(k,x)}\right)\right)\\
&\to& H^0\left(Y, \sigma^*\sO_X(-kr_0K_X)\cdot\left(I_{(k,x)}\cdot\sO_Y\right)\right)\\
&\hookrightarrow& H^0\left(Y, \sigma^*\sO_X(-kr_0K_X)(\lfloor -xF\rfloor)\right).
\end{eqnarray*}
Thus we get $H^0(X, \sO_X(-kr_0K_X)\cdot I_{(k,x)})\subset \sF^xV_k$.
\end{proof}

Take any $e_+$, $e_-\in\Z$ with $e_+>r_0\cdot\tau(F)$ and $e_-<0$. 
Let $r_1\in\Z_{>0}$ be a sufficiently big positive integer 
as in \cite[Proposition 4.3 (4)]{fjt2}. For any $r\geq r_1$, we set the flag ideal 
$\sI_r\subset\sO_{X\times\A^1_t}$ (in the sense of \cite[Definition 3.1]{odk}) defined by 
\[
\sI_r:=I_{(r, re_+)}+I_{(r, re_+-1)}t^1+\cdots+I_{(r, re_-+1)}t^{r(e_+-e_-)-1}+
(t^{r(e_+-e_-)}).
\]
Let $\Pi_r\colon\sX^r\to X\times\A^1$ be the blowup along $\sI_r$, let 
$E^r\subset\sX^r$ be the Cartier divisor defined by 
$\sO_{\sX^r}(-E^r)=\sI_r\cdot\sO_{\sX^r}$, and set 
$\sL^r:=\Pi^*_r(-K_{X\times\A^1})-(1/(rr_0))E^r$. 
Then $(\sX^r, \sL^r)/\A^1$ is a semiample test configuration for $(X, -K_X)$ by 
\cite[Lemma 4.6]{fjt2}.

\begin{claim}\label{lambda_claim}
We have
\[
\lim_{r\to\infty}\lambda_{\max}(\sX^r, \sL^r)=-\frac{e_+}{r_0}+\tau(F).
\]
\end{claim}

\begin{proof}[Proof of Claim \ref{lambda_claim}]
Since $rr_0\sL^r$ is Cartier, we can consider the filtration 
$\sF_{(\sX^r, rr_0\sL^r)}$ of $V_\bullet^{(r)}:=\bigoplus_{k\in\Z_{\geq 0}}V_{kr}$ 
as in Proposition \ref{tcfilt_prop}. For any $k\in\Z_{\geq 0}$, we set 
$J_{(k;r,kre_-)}:=\sO_X$ and 
\[
\sI_r^k=:J_{(k;r, kre_+)}+J_{(k;r,kre_+-1)}t^1+\cdots
+J_{(k;r,kre_-+1)}t^{kr(e_+-e_-)-1}+(t^{kr(e_+-e_-)}). 
\]
Then, by Proposition \ref{tcfilt_prop} \eqref{tcfilt_prop1}, 
$\sF_{(\sX^r, rr_0\sL^r)}^xV_k^{(r)}$ is equal to 
\[\begin{cases}
V_{kr} & \text{if }x\leq -kr(e_+-e_-),\\
H^0(X, \sO_X(-krr_0K_X)\cdot J_{(k;r,kre_++\lceil x\rceil)}) & \text{if }-kr(e_+-e_-)
<x\leq 0,\\
0 & \text{if }x>0
\end{cases}\]
for $k\gg 0$. By \cite[Lemma 4.5]{fjt2}, for $j\in[kre_-,kre_+]\cap\Z$, 
\[
H^0(X, \sO_X(-krr_0K_X)\cdot J_{(k;r,j)})=0
\]
holds if and only if $J_{(k;r,j)}=0$ holds. 
On the other hand, 
\begin{eqnarray*}
&&\max\{j\in[kre_-,kre_+]\cap\Z\,|\, J_{(k;r,j)}\neq 0\}\\
&=&k\cdot\max\{j\in[re_-,re_+]\cap\Z\,|\, I_{(r,j)}\neq 0\}\\
&=&k\cdot\max\{j\in[re_-,re_+]\cap\Z\,|\, \sF^jV_r\neq 0\}
\end{eqnarray*}
by Claim \ref{sat_claim}. This implies that 
\begin{eqnarray*}
&&\sup\{x\in\R\,|\,\sF^x_{(\sX^r, rr_0\sL^r)}V_k^{(r)}\neq 0\}\\
&=&k\left(-re_++\max\{j\in[re_-,re_+]\cap\Z\,|\, \sF^jV_r\neq 0\}\right).
\end{eqnarray*}
By Proposition \ref{tcfilt_prop} \eqref{tcfilt_prop2}, we get the assertion.
\end{proof}

Let $\nu\colon\sX^{r,\nu}\to\sX^r$ be the normalization. We know
\[
\Ding(\sX^{r,\nu},\nu^*\sL^r)\geq \delta\cdot J^{\NA}(\sX^{r,\nu},\nu^*\sL^r)
\]
by the assumption. Note that 
\[
\lct(\sX^{r,\nu}, D_{(\sX^{r,\nu},\nu^*\sL^r)};\sX^{r,\nu}_0)
=\lct(X\times\A^1, \sI_r^{\cdot 1/(rr_0)};(t))
\]
by the proof of \cite[Proposition 3.5]{fjt2}. Thus 
\[
\lct(X\times\A^1, \sI_r^{\cdot 1/(rr_0)};(t))\geq d_{r,\delta}, 
\]
where
\[
d_{r,\delta}:=1+\frac{(\bar{\sL}^{r\cdot n+1})}{(n+1)((-K_X)^{\cdot n})}
+\delta\cdot\left(\lambda_{\max}(\sX^r, \sL^r)
-\frac{(\bar{\sL}^{r\cdot n+1})}{(n+1)((-K_X)^{\cdot n})}\right).
\]
Set $F_0:=F$, let $\{F_i\}_{i\in I}$ be the set of $\sigma$-exceptional prime divisors on 
$Y$ and set $\hat{I}:=I\cup\{0\}$. Then the pair 
\[
\biggl(Y\times\A^1, \prod_{i\in\hat{I}}\sO_{Y\times\A^1}(-F_i\times\A^1
)^{\cdot(-A_X(F_i)+1)}\cdot(\sI_r\cdot\sO_{Y\times\A^1})^{\cdot 1/(rr_0)}
\cdot(t)^{\cdot d_{r,\delta}}\biggr)
\]
is sub log canonical for any $r\geq r_1$. We know that 
\begin{eqnarray*}
\sI_r\cdot\sO_{Y\times\A^1}&\subset&\sO_Y(-re_+F)+\sO_Y(-(re_+-1)F)t^1+\cdots\\
&&\cdots+\sO_Y(-F)t^{re_+-1}+(t^{re_+})=\left(\sO_Y(-F)+(t)\right)^{re_+}.
\end{eqnarray*}

\begin{claim}\label{ddelta_claim}
The limit $d_{\infty,\delta}:=\lim_{r\to\infty}d_{r,\delta}$ exists and is equal to 
\[
1-\frac{e_+}{r_0}+\delta\cdot\tau(F)+\frac{1-\delta}{((-K_X)^{\cdot n})}
\int_0^\infty\vol_X(-K_X-xF)dx.
\]
\end{claim}

\begin{proof}[Proof of Claim \ref{ddelta_claim}]
By \cite[Lemma 4.7]{fjt2}, we have 
\begin{eqnarray*}
\lim_{r\to\infty}\frac{(\bar{\sL}^{r\cdot n+1})}{(n+1)((-K_X)^{\cdot n})}
&=&-\frac{e_+-e_-}{r_0}\\
&&+\frac{1}{r_0^{n+1}((-K_X)^{\cdot n})}
\int_{e_-}^{e_+}\vol(\sF V_\bullet^x)dx.
\end{eqnarray*}
Note that $\vol(\sF V_\bullet^x)$ is equal to 
\[\begin{cases}
r_0^n\cdot\vol_X(-K_X-\frac{x}{r_0}F) & \text{if }x\in\R_{\geq 0}, \\
r_0^n\cdot((-K_X)^{\cdot n}) & \text{otherwise}.
\end{cases}\]
Thus Claim \ref{ddelta_claim} follows from Claim \ref{lambda_claim}.
\end{proof}

By Claim \ref{ddelta_claim}, the pair 
\[
\biggl(Y\times\A^1, \prod_{i\in\hat{I}}\sO_{Y\times\A^1}(-F_i\times\A^1
)^{\cdot(-A_X(F_i)+1)}\cdot(\sO_Y(-F)+(t))^{\cdot e_+/r_0}
\cdot(t)^{\cdot d_{\infty,\delta}}\biggr)
\]
is sub log canonical. 
Take the blowup $\sY\to Y\times\A^1$ of $Y\times\A^1$ along $F\times\{0\}$, 
and let $E_F$ be the exceptional divisor of the blowup. Then we have 
\begin{eqnarray*}
-1&\leq&\ord_{E_F}(K_{\sY/Y\times\A^1})-(-A_X(F)+1)-\frac{e_+}{r_0}-d_{\infty,\delta}\\
&=&A_X(F)-\frac{e_+}{r_0}-d_{\infty,\delta}.
\end{eqnarray*}
By Claim \ref{ddelta_claim}, we get 
\[
\delta\cdot\tau(F)+\frac{1-\delta}{((-K_X)^{\cdot n})}
\int_0^\infty\vol_X(-K_X-xF)dx\leq A_X(F).
\]
This immediately implies that $\beta(F)\geq \delta\cdot j(F)$.
\end{proof}

\section{Special test configurations and dreamy prime divisors}\label{vK_section}

In this section, we prove Theorems \ref{mainthm} and \ref{mainK_thm}. 
More precisely, we prove the following: 

\begin{thm}\label{vK_thm}
Let $X$ be a $\Q$-Fano variety and let $(\sX, -K_{\sX/\A^1})/\A^1$ be a normal, 
ample test configuration for $(X, -K_X)$ with $\sX_0$ irreducible and reduced. 
Then the divisorial valuation $v_{\sX_0}$ on $X$ defined in Proposition \ref{discrep_prop}
is dreamy over $X$ and we have the equalities 
$\DF(\sX, -K_{\sX/\A^1})=\beta(v_{\sX_0})/\vol_X(-K_X)$ and 
$J^{\NA}(\sX, -K_{\sX/\A^1})=j(v_{\sX_0})/\vol_X(-K_X)$. 
\end{thm}

\begin{thm}\label{K_thm}
Let $X$ be a $\Q$-Fano variety and $F$ a dreamy prime divisor over $X$. 
Then there exists a semiample test configuration $(\sX, \sL)/\A^1$ for $(X, -K_X)$ 
such that the equalities $\DF(\sX, \sL)=\beta(F)/\vol_X(-K_X)$ and 
$J^{\NA}(\sX, \sL)=j(F)/\vol_X(-K_X)$ 
hold. 
\end{thm}

\begin{remark}\label{vK_rmk}
From Theorem \ref{vK_thm}, it seems to be more natural to consider the 
invariants $\beta(v)/\vol_X(-K_X)$ and $j(v)/\vol_X(-K_X)$ rather than 
$\beta(v)$ and $j(v)$ (see also \cite[Definition 3.4]{li2}). The reason why the author 
uses the invariants $\beta(v)$ and $j(v)$ is because the invariant $\beta(v)$ is a 
modification of the invariants in \cite[Proposition 3.2]{slope}, 
\cite[Definition 1.1]{fjt1}, \cite[\S 3]{fjtdP} and \cite[Theorem 4.10]{fjt2} 
introduced by the author. 
\end{remark}

\begin{proof}[Proof of Theorem \ref{mainthm}]
This is an immediate consequence of Corollary \ref{stc_cor}, 
Theorems \ref{Kv_thm} and \ref{vK_thm}.
\end{proof}

\begin{proof}[Proof of Theorem \ref{mainK_thm}]
This is an immediate consequence of Theorems \ref{vK_thm} and \ref{K_thm}. 
\end{proof}

\begin{proof}[Proof of Theorem \ref{vK_thm}]
Pick $r_0\in\Z_{>0}$ with $-r_0K_{\sX/\A^1}$ Cartier. Let 
\[\xymatrix{
& \bar{\sZ}  \ar[dl]_\Pi \ar[dr]^\Theta & \\
X\times\pr^1 & & \bar{\sX}
}\]
be the normalization of the graph. Let $\sZ_0=\sum_{i\in I}m_iE_i+\hat{X}_0+\hat{\sX}_0$ 
be the irreducible decomposition, 
where $\hat{X}_0$ is the strict transform of $X\times\{0\}$ and 
$\hat{\sX}_0$ is the strict transform of $\sX_0$. We set 
$B:=\Theta^*(-K_{\bar{\sX}/\pr^1})-\Pi^*p_1^*(-K_X)$ supported on $\sZ_0$. Note that 
$-\ord_{\hat{\sX}_0}B=\ord_{\hat{\sX}_0}(K_{\sZ/X\times\A^1})$. 
Let $V_\bullet$ be the complete graded linear series of $-r_0K_X$ and let us consider 
the filtration $\sF:=\sF_{(\sZ, \Theta^*(-r_0K_{\sX/\A^1}))}$ of $V_\bullet$ as 
in Proposition \ref{tcfilt_prop}.

\begin{claim}[{cf.\ \cite[Lemma 5.16]{BHJ} and \cite[Lemma 6.6]{li2}}]\label{rest_claim}
Pick $k\in\Z_{\geq 0}$ and $x\in\R$. Then we have 
\[
\sF^xV_k=\{f\in V_k\,|\, v_{\sX_0}(f)\geq kr_0A_X(v_{\sX_0})+x\}.
\]
\end{claim}

\begin{proof}[Proof of Claim \ref{rest_claim}]
Take any $f\in V_k\setminus\{0\}$. Let $D\in|-kr_0K_X|$ be the effective divisor 
which corresponds to $f$. Then $f\in\sF^xV_k$ holds if and only if 
$\Pi^*p_1^*D+kr_0B\geq \lceil x\rceil\sZ_0$ holds by Proposition \ref{tcfilt_prop} 
\eqref{tcfilt_prop1}. Since 
\[
\Pi^*p_1^*D+kr_0B- \lceil x\rceil\sZ_0\sim_\Q
\Theta^*(-kr_0K_{\bar{\sX}/\pr^1}- \lceil x\rceil\sX_0), 
\]
the condition $\Pi^*p_1^*D+kr_0B- \lceil x\rceil\sZ_0\geq 0$ 
is equivalent to the condition 
$\Theta_*(\Pi^*p_1^*D+kr_0B- \lceil x\rceil\sZ_0)\geq 0$. 
This condition is equivalent to the condition 
\[
\ord_{\hat{\sX}_0}(\Pi^*p_1^*D)+kr_0\cdot\ord_{\hat{\sX}_0}B- \lceil x\rceil\geq 0.
\]
Note that $-\ord_{\hat{\sX}_0}B=A_X(v_{\sX_0})$ by Proposition \ref{discrep_prop}.
\end{proof}

By Claim \ref{rest_claim}, $\sF^xV_k$ is equal to 
\[\begin{cases}
H^0(X, -kr_0K_X-(kr_0A_X(v_{\sX_0})+x)v_{\sX_0}) & \text{if }x\geq 
-kr_0A_X(v_{\sX_0}),\\
V_k & \text{otherwise.}
\end{cases}\]
We note that the $\Bbbk$-algebra $\bigoplus_{k\in\Z_{\geq 0}, j\in\Z}\sF^jV_k$ is 
finitely generated by Proposition \ref{tcfilt_prop} \eqref{tcfilt_prop1}. 
Thus the $\Bbbk$-algebra
\[
\bigoplus_{k,j\in\Z_{\geq 0}}H^0(X, -kr_0K_X-jv_{\sX_0})
\]
is finitely generated (see \cite[Lemma 4.8]{ELMNP}). 
In particular, the divisorial valuation $v_{\sX_0}$ is dreamy. 
By Proposition \ref{tcfilt_prop}, 
\[
\lambda_{\max}(\sX, -K_{\sX/\A^1})=\tau(v_{\sX_0})-A_X(v_{\sX_0}), 
\]
and 
\[
\lambda_{\min}^{(k)}:=\inf\{x\in\R\,|\,\sF^xV_k\neq V_k\}
\]
satisfies that 
\[
\lim_{k\to\infty}\frac{\lambda_{\min}^{(k)}}{kr_0}=-A_X(v_{\sX_0}). 
\]
Then the quantity $w(k)$ in Proposition \ref{tcfilt_prop} is equal to 
\[
\int_0^{kr_0\tau(v_{\sX_0})}\dim H^0(X, -kr_0K_X-xv_{\sX_0})dx
-kr_0A_X(v_{\sX_0})\dim V_k. 
\]
Thus, by Proposition \ref{tcfilt_prop} 
\eqref{tcfilt_prop22}, $((-K_{\bar{\sX}/\pr^1})^{\cdot n+1})$ is equal to
\begin{eqnarray*}
&&\lim_{k\to\infty}\frac{w(k)}
{(kr_0)^{n+1}/(n+1)!}\\
&=&(n+1)\left(\int_0^{\tau(v_{\sX_0})}\vol_X(-K_X-xv_{\sX_0})dx
-A_X(v_{\sX_0})((-K_X)^{\cdot n})\right)\\
&=&-(n+1)\beta(v_{\sX_0}).
\end{eqnarray*}
Hence 
\begin{eqnarray*}
\DF(\sX, -K_{\sX/\A^1})&=&\frac{\beta(v_{\sX_0})}{((-K_X)^{\cdot n})}, \\
J^{\NA}(\sX, -K_{\sX/\A^1})&=&\tau(v_{\sX_0})-A_X(v_{\sX_0})+
\frac{\beta(v_{\sX_0})}{((-K_X)^{\cdot n})} \\
&=&\frac{j(v_{\sX_0})}{((-K_X)^{\cdot n})}.
\end{eqnarray*}
Therefore we get the assertion. 
\end{proof}

\begin{proof}[Proof of Theorem \ref{K_thm}]
The\, proof\, is\, essentially\, same\, as\, the\, one\, in\, \cite[Theorem 1.1]{fjtdP}. 
Take any projective log resolution $\sigma\colon Y\to X$ with $F$ a smooth divisor 
on $Y$. Fix a sufficiently divisible positive integer $r_0\in\Z_{>0}$ such that 
$-r_0K_X$ is Cartier and the $\Bbbk$-algebra 
\[
\bigoplus_{k,j\in\Z_{\geq 0}}H^0(Y, \sigma^*(-kr_0K_X)-jF)
\]
is generated by 
\[
\bigoplus_{j\in\Z_{\geq 0}}H^0(Y, \sigma^*(-r_0K_X)-jF).
\]
We set $H:=\sigma^*(-K_X)$ and $V_{k,j}:=H^0(Y, kr_0H-jF)$ for simplicity. 
We set 
\[
I_j:=\Image\left(V_{1,j}\otimes_\Bbbk\sO_X(r_0K_X)\to\sO_X\right)
\]
for $j\in[0,r_0\tau(F)]\cap\Z$, and set 
\[
\sI:=I_{r_0\tau(F)}+I_{r_0\tau(F)-1}t^1+\cdots+I_1t^{r_0\tau(F)-1}+(t^{r_0\tau(F)})
\subset\sO_{X\times\A^1_t}.
\]
Let $\Pi\colon\sX\to X\times\A^1$ be the blowup along $\sI$, let $E\subset\sX$ be 
the Cartier divisor defined by $\sO_\sX(-E)=\sI\cdot\sO_\sX$, and let 
$\sL:=\Pi^*p_1^*(-K_X)-(1/r_0)E$. 

\begin{claim}[{\cite{fjtdP}}]\label{fg_claim}
\begin{enumerate}
\renewcommand{\theenumi}{\arabic{enumi}}
\renewcommand{\labelenumi}{(\theenumi)}
\item\label{fg_claim1}
For any $k\in\Z_{>0}$, set $J_{(k,0)}:=\sO_X$ and 
\[
\sI^k=:J_{(k,kr_0\tau(F))}+J_{(k,kr_0\tau(F)-1)}t^1+\cdots+J_{(k,1)}t^{kr_0\tau(F)-1}
+(t^{kr_0\tau(F)}).
\]
Then 
\[
V_{k,j}=H^0(X, \sO_X(-kr_0K_X)\cdot J_{(k,j)})
\]
holds for any $j\in[0,kr_0\tau(F)]\cap\Z$. 
\item\label{fg_claim2}
$(\sX, \sL)/\A^1$ is a semiample test configuration for $(X, -K_X)$.
\end{enumerate}
\end{claim}

\begin{proof}[Proof of Claim \ref{fg_claim}]
\eqref{fg_claim1} follows from the proof of \cite[Lemma 3.6]{fjtdP}, \eqref{fg_claim2} 
follows from the proof of \cite[Lemma 3.7]{fjtdP}.
\end{proof}

Let $V_\bullet$ be the complete graded linear series of $-r_0K_X$ and let us consider 
the filtration $\sF:=\sF_{(\sX, r_0\sL)}$ of $V_\bullet$ as in Proposition \ref{tcfilt_prop}. 
By Proposition \ref{tcfilt_prop} \eqref{tcfilt_prop1} and Claim \ref{fg_claim} 
\eqref{fg_claim1}, for $k\gg 0$, we have 
\[
\sF^xV_k=\begin{cases}
V_k & \text{if }x\leq -kr_0\tau(F), \\
V_{k,kr_0\tau(F)+\lceil x\rceil} & \text{if }-kr_0\tau(F)<x\leq 0,\\
0 & \text{if }x>0.
\end{cases}
\]
Thus $w(k)$ in Proposition \ref{tcfilt_prop} is equal to 
$f(k)-kr_0\tau(F)\dim V_k$, where 
\[
f(k):=\sum_{l=1}^{kr_0\tau(F)}\dim V_{k,l}.
\]
By Proposition \ref{tcfilt_prop} \eqref{tcfilt_prop2}, $f(k)$ is a polynomial function of 
degree at most $n+1$ for $k\gg 0$. Let us write 
$f(k)=f_{n+1}k^{n+1}+f_nk^n+O(k^{n-1})$. By the asymptotic Riemann-Roch Theorem, 
we know that 
\[
\dim V_k=\frac{r_0^n((-K_X)^{\cdot n})}{n!}k^n+\frac{r_0^{n-1}((-K_X)^{\cdot n})}
{2\cdot((n-1)!)}k^{n-1}+O(k^{n-2})
\]
for $k\gg 0$. Thus, by Proposition \ref{tcfilt_prop} \eqref{tcfilt_prop23} and 
\eqref{tcfilt_prop24}, we get
\[
\DF(\sX, \sL)=\frac{n!}{((-r_0K_X)^{\cdot n})}\left(\frac{n}{r_0}f_{n+1}-2f_n\right).
\]

By \cite[Theorem 4.2]{KKL}, there exist a sequence of rational numbers 
\[
0=\tau_0<\tau_1<\cdots<\tau_m=\tau(F)
\]
and pairwise distinct birational contraction maps 
$\varphi_j\colon Y\dashrightarrow Y_j$ with $Y_j$ normal and projective for 
$1\leq j\leq m$ such that the map $\varphi_j$ is a semiample model 
of $H-xF$ for any $x\in[\tau_{j-1}, \tau_j]$, and the ample model of $H-xF$ for 
any $x\in(\tau_{j-1}, \tau_j)$. See \cite[Definition 2.3]{KKL} for the definitions of 
semiample and ample models. By \cite[Remark 2.4 (i)]{KKL}, we have 
\[
f(k)=\sum_{j=1}^m\sum_{l=kr_0\tau_{j-1}+1}^{kr_0\tau_j}\dim 
H^0(Y_j, kr_0H_j-lF_j), 
\]
where $H_j$, $F_j$ is the strict transform of $H$, $F$ on $Y_j$, respectively. 
Hence, by \cite[Proposition 4.1]{fjt1}, we have 
\begin{eqnarray*}
f(k)&=&\sum_{j=1}^m\biggl(\frac{(kr_0)^{n+1}}{n!}\int_{\tau_{j-1}}^{\tau_j}
((H_j-xF_j)^{\cdot n})dx\\ 
&-&\frac{(kr_0)^n}{2\cdot (n-1)!}\int_{\tau_{j-1}}^{\tau_j}
((H_j-xF_j)^{\cdot n-1}\cdot (K_{Y_j}+F_j))dx\biggr)+O(k^{n-1}).
\end{eqnarray*}

\begin{claim}\label{exc_claim}
For any $\sigma$-exceptional prime divisor $F'$ on $Y$ with $F'\neq F$ and for 
any $1\leq j\leq m$, the divisor $F'$ is $\varphi_j$-exceptional. In particular, 
$H_j+K_{Y_j}+F_j=A_X(F)\cdot F_j$ holds for any $1\leq j\leq m$.
\end{claim}

\begin{proof}[Proof of Claim \ref{exc_claim}]
Fix any $\tau\in(\tau_{j-1}, \tau_j)\cap\Q$. Let 
\[\xymatrix{
& Z  \ar[dl]_\pi \ar[dr]^\theta & \\
Y & & Y_j
}\]
be a common resolution of $\varphi_j$. Then there exists an ample $\Q$-divisor 
$A_j$ on $Y_j$ such that the linear system $|\pi^*(k(H-\tau F))|$ is equal to 
$|\theta^*kA_j|+G_k$ with $G_k$ effective and fixed for any sufficiently divisible 
$k\in\Z_{>0}$. In particular, $G_k$ is $\theta$-exceptional. Let $F'_Z\subset Z$ be 
the strict transform of $F'$. Assume that $F'_Z$ is not $\theta$-exceptional. 
Since $F'_Z$ is $(\sigma\circ\pi)$-exceptional, $F'_Z$ is covered by a family of curves 
$\{C_t\}_{t\in T}$ with $(\pi^*(H-\tau F)\cdot C_t)\leq 0$. However, a general 
$C_t$ satisfies that $(\theta^*A_j\cdot C_t)=(A_j\cdot \theta_*C_t)>0$ and 
$(G_k\cdot C_t)\geq 0$. This implies that $(\pi^*(H-\tau F)\cdot C_t)>0$, 
a contradiction. 
The equality $H_j+K_{Y_j}+F_j=A_X(F)\cdot F_j$ follows immediately from the fact 
that all of those $F'$ are $\varphi_j$-exceptional. 
\end{proof}

By Claim \ref{exc_claim} and \cite[Remark 2.4 (i)]{KKL}, we get 
\begin{eqnarray*}
&&\DF(\sX, \sL)\\
 &=& \frac{n}{((-K_X)^{\cdot n})}\sum_{j=1}^m
\int_{\tau_{j-1}}^{\tau_j}(A_X(F)-x)\left((H_j-xF_j)^{\cdot n-1}\cdot F_j\right)dx\\
&=&\frac{1}{((-K_X)^{\cdot n})}\sum_{j=1}^m
\biggl([(x-A_X(F))((H_j-xF_j)^{\cdot n})]^{\tau_j}_{\tau_{j-1}}\\
&&
-\int_{\tau_{j-1}}^{\tau_j}((H_j-xF_j)^{\cdot n})dx\biggr)=\frac{\beta(F)}{((-K_X)^{\cdot n})}.
\end{eqnarray*}
We can also check that 
\begin{eqnarray*}
\lambda_{\max}(\sX, \sL)&=&0, \\
(\bar{\sL}^{\cdot n+1})&=&(n+1)!\cdot\frac{f_{n+1}}{r_0^{n+1}}
-(n+1)\tau(F)((-K_X)^{\cdot n}),\\
f_{n+1} &=& \sum_{j=1}^m\frac{r_0^{n+1}}{n!}\int_{\tau_{j-1}}^{\tau_j}
((H_j-xF_j)^{\cdot n})dx.
\end{eqnarray*}
Thus we get the equality $J^{\NA}(\sX, \sL)=j(F)/((-K_X)^{\cdot n})$.
\end{proof}

\section{Log Fano pairs}\label{log_section}

We can naturally and easily generalize Theorem \ref{mainthm} for log Fano pairs. 
We omit the proofs of many results in this section since the results are direct 
generalizations of the results in \S \ref{prelim_section}--\S \ref{vK_section}. 
We remark that some of the results in this section had already been known 
(e.g., \cite[\S 4]{LL}, \cite[\S 4]{LX2}) after the author 
uploaded the first version of this article on arXiv.

\begin{definition}[{see Definitions \ref{intro_dfn1} and \ref{intro_dfn2}}]\label{log_intro_dfn}
Let $(X, \Delta)$ be a \emph{log Fano pair}, that is, 
$(X, \Delta)$ is a projective klt pair with $\Delta$ effective 
$\Q$-divisor and $-(K_X+\Delta)$ ample $\Q$-Cartier. 
Given a prime divisor $F$ over $X$, the definitions of the quantities $A_{(X, \Delta)}(F)$, 
$\tau_{(X, \Delta)}(F)$, $\beta_{(X, \Delta)}(F)$ and $j_{(X, \Delta)}(F)$ are exactly the 
same as in Definition \ref{intro_dfn2}, replacing $K_X$ with $K_X+\Delta$ in all 
formulas. The same is true for the definition of $F$ to be dreamy over $(X, \Delta)$. 
\end{definition}

\begin{definition}[{see Definition \ref{triv_dfn}}]\label{log_triv_dfn}
Let $(X, \Delta)$ be a log Fano pair and let 
$(\sX, \sL)/\A^1$ be a normal test configuration for $(X, -(K_X+\Delta))$. 
\begin{enumerate}
\renewcommand{\theenumi}{\arabic{enumi}}
\renewcommand{\labelenumi}{(\theenumi)}
\item\label{log_triv_dfn1}
Let $\Delta_{\sX}$ (resp.\ $\Delta_{\bar{\sX}}$) be the 
$\Q$-divisor which is the closure of $\Delta\times(\A^1\setminus\{0\})$ on 
$\sX$ (resp.\ on $\bar{\sX}$).
\item\label{log_triv_dfn2}
A test configuration $(\sX, \sL)/\A^1$ is said to be a \emph{product-type} 
test configuration for $((X, \Delta), -(K_X+\Delta))$ if 
$(\sX, \Delta_{\sX})$ is isomorphic to $(X\times\A^1, \Delta\times\A^1)$. 
\item\label{log_triv_dfn3}
A test configuration 
$(\sX, \sL)/\A^1$ is said to be a \emph{special test configuration} for 
$((X, \Delta), -(K_X+\Delta))$ 
if it is a normal, ample test configuration for $(X, -(K_X+\Delta))$, 
$\sL=-(K_{\sX/\A^1}+\Delta_{\sX})$ and the pair 
$(\sX,\Delta_{\sX}+\sX_0)$ is plt. 
\end{enumerate}
\end{definition}

We recall the notions of Ding invariants and Donaldson-Futaki invariants. 
See also \cite[Definitions 3.17 and 7.25]{BHJ} and \cite[Definitions 2.2 and 2.3]{LL}. 

\begin{definition}[{see Definition \ref{DF_dfn}}]\label{log_DF_dfn}
Let $(X, \Delta)$ be a log Fano pair of dimension $n$ and $(\sX, \sL)/\A^1$ be a 
normal test configuration for $(X, -(K_X+\Delta))$. 
\begin{enumerate}
\renewcommand{\theenumi}{\arabic{enumi}}
\renewcommand{\labelenumi}{(\theenumi)}
\item\label{log_DF_dfn3}
Let $D_{((\sX,\Delta_{\sX}), \sL)}$ be the $\Q$-divisor on $\bar{\sX}$ defined by 
\begin{itemize}
\item
$\Supp D_{((\sX,\Delta_{\sX}), \sL)}\subset\sX_0$, 
\item
$D_{((\sX,\Delta_{\sX}), \sL)}\sim_\Q 
-(K_{\bar{\sX}/\pr^1}+\Delta_{\bar{\sX}})-\bar{\sL}$.
\end{itemize}
The \emph{Ding invariant} $\Ding_\Delta(\sX, \sL)$ of 
$(\sX, \sL)/\A^1$ is defined by the 
following: 
\begin{eqnarray*}
\Ding_\Delta(\sX, \sL)&:=&-\frac{(\bar{\sL}^{\cdot n+1})}
{(n+1)((-(K_X+\Delta))^{\cdot n})}-1\\
&&+\lct(\sX, \Delta_{\sX}+D_{((\sX, \Delta_{\sX}), \sL)}; \sX_0),
\end{eqnarray*}
where
\begin{eqnarray*}
&&\lct(\sX, \Delta_{\sX}+D_{((\sX, \Delta_{\sX}), \sL)}; \sX_0)\\
&:=&\max\{c\in\R\,|\,(\sX, \Delta_{\sX}+D_{((\sX, \Delta_{\sX}), \sL)}+c\sX_0)\text{: 
 sub log canonical}\}.
\end{eqnarray*}
\item\label{log_DF_dfn4}
The \emph{Donaldson-Futaki invariant} $\DF_\Delta(\sX, \sL)$ 
of the test configuration $(\sX, \sL)/\A^1$ is 
defined by the following: 
\begin{eqnarray*}
\DF_\Delta(\sX, \sL)&:=&\frac{n}{n+1}\cdot\frac{(\bar{\sL}^{\cdot n+1})}
{((-(K_X+\Delta))^{\cdot n})}\\
&&+\frac{(\bar{\sL}^{\cdot n}\cdot (K_{\bar{\sX}/\pr^1}+\Delta_{\bar{\sX}}))}
{((-(K_X+\Delta))^{\cdot n})}.
\end{eqnarray*}
\end{enumerate}
\end{definition}

We note that the analogue of Proposition \ref{pull_prop} holds.

\begin{definition}[{see Definition \ref{K_dfn}}]\label{log_K_dfn}
Let $(X, \Delta)$ be a log Fano pair. 
The definition of uniform K-stability, K-stability, K-polystability, K-semistability, 
uniform Ding stability, Ding stability, Ding polystability, Ding semistability of $(X, \Delta)$ 
is the same as in Definition \ref{K_dfn}. For example, $(X, \Delta)$ is K-polystable if 
$\DF(\sX, \sL)\geq 0$ for any normal, ample test configuration $(\sX, \sL)/\A^1$ 
for $(X, -(K_X+\Delta))$ and equality holds only if $(\sX, \sL)/\A^1$ is a 
product-type for $((X, \Delta), -(K_X+\Delta))$. 
\end{definition}

The main theorem in this section is the following: 

\begin{thm}\label{log_mainthm}
Let $(X, \Delta)$ be a log Fano pair. For any $\delta\in[0, 1)$, the following are 
equivalent: 
\begin{enumerate}
\renewcommand{\theenumi}{\roman{enumi}}
\renewcommand{\labelenumi}{(\theenumi)}
\item\label{log_mainthm_1}
$\DF_\Delta(\sX, \sL)\geq \delta\cdot J^{\NA}(\sX, \sL)$ holds for any normal, ample 
test configuration $(\sX, \sL)/\A^1$ for $(X, -(K_X+\Delta))$.
\item\label{log_mainthm_2}
$\Ding_\Delta(\sX, \sL)\geq \delta\cdot J^{\NA}(\sX, \sL)$ holds for any normal, ample 
test configuration $(\sX, \sL)/\A^1$ for $(X, -(K_X+\Delta))$.
\item\label{log_mainthm_3}
$\beta_{(X, \Delta)}(F)\geq \delta\cdot j_{(X, \Delta)}(F)$ holds 
for any prime divisor $F$ over $X$.
\item\label{log_mainthm_4}
$\beta_{(X, \Delta)}(F)\geq \delta\cdot j_{(X, \Delta)}(F)$ holds 
for any dreamy prime divisor $F$ over $(X, \Delta)$.
\end{enumerate}
\end{thm}

\begin{example}[{\cite[Theorem 3]{liR}}]\label{logP1_ex}
We consider log Fano pairs of dimension one. 
Let us consider a pair $(\pr^1, \Delta)$ with $\Delta=\sum_{i=1}^ma_ip_i$ 
$(m\geq 1$, $a_i\in\Q_{>0}$, $p_1,\dots,p_m$ are distinct points$)$, 
$1>a_1\geq\cdots\geq a_m$ and $2-\sum_{i=1}^ma_i>0$. 
We remark that any prime divisor over $\pr^1$ is a point on $\pr^1$. Thus we can 
check from Theorem \ref{log_mainthm} that 
$(\pr^1, \Delta)$ is uniformly K-stable 
(resp.\ K-semistable) if and only if 
$\sum_{i=2}^ma_i>a_1$ (resp.\ $\sum_{i=2}^ma_i\geq a_1$) holds. 
\end{example}

The following three theorems are logarithmic versions of the results in \cite{LX} and 
\cite{BBJ}. 

\begin{thm}[{see Theorem \ref{LX2_thm}}]\label{log_LX2_thm}
Let $(X, \Delta)$ be a  log Fano pair and $(\sX, \sL)/\A^1$ be a normal, ample 
test configuration for $(X, -(K_X+\Delta))$. 
Then\, there\, exist\, $d\in\Z_{>0}$,\, a\, projective\, birational\, 
$\G_m$-equivariant\, morphism\, $\pi\colon \sX^{\LC}\to\sX^{(d)}$,\, and\, 
a\, normal,\, ample\, test\, configuration\, 
$(\sX^{\LC}, \sL^{\LC})/\A^1$ for $(X, -(K_X+\Delta))$ 
such that: 
\begin{enumerate}
\renewcommand{\theenumi}{\arabic{enumi}}
\renewcommand{\labelenumi}{(\theenumi)}
\item\label{log_LX2_thm1}
$(\sX^{\LC}, \Delta_{\sX^{\LC}}+\sX_0^{\LC})$ is log canonical. 
\item\label{log_LX2_thm2}
$\DF_\Delta(\sX^{\LC}, \sL^{\LC})\leq d\cdot\DF_\Delta(\sX, \sL)$ holds. 
Moreover, 
equality holds if and only if $(\sX, \Delta_{\sX}+\sX_0)$ is log canonical and 
$-(K_{\sX^{(d)}}+\Delta_{\sX^{(d)}})\sim_\Q 
\psi_d^*\sL$; in which case $\sL^{\LC}\sim_\Q\pi^*\psi_d^*\sL$ and $\pi$
is an isomorphism. 
\item\label{log_LX2_thm3}
$\Ding_\Delta(\sX^{\LC}, \sL^{\LC})\leq d\cdot\Ding_\Delta(\sX, \sL)$ holds. 
Moreover, 
equality holds if and only if $(\sX^{(d)}, \Delta_{\sX^{(d)}}+\sX_0^{(d)})$ 
is log canonical and 
$-(K_{\sX^{(d)}}+\Delta_{\sX^{(d)}})\sim_\Q\psi_d^*\sL$; 
in which case $\sL^{\LC}\sim_\Q\pi^*\psi_d^*\sL$ and $\pi$
is an isomorphism. 
\item\label{log_LX2_thm4}
For any $\delta\in[0,1]$, we have the inequality 
\[
d(\Ding_\Delta(\sX, \sL)-\delta\cdot J^{\NA}(\sX, \sL))\geq
\Ding_\Delta(\sX^{\LC}, \sL^{\LC})-\delta\cdot J^{\NA}(\sX^{\LC}, \sL^{\LC}).
\]
\end{enumerate}
\end{thm}

\begin{proof}
We follow the proof of Proposition \ref{LX2_thm}. 
We take the log canonical modification $\pi\colon \sX^{\LC}\to (\sX^{(d)}, 
\Delta_{\sX^{(d)}}+\sX_0^{(d)})$ of the pair 
$(\sX^{(d)}, \Delta_{\sX^{(d)}}+\sX_0^{(d)})$ for some $d\in\Z_{>0}$.
\end{proof}

\begin{thm}[{see Theorem \ref{LX3_thm}}]\label{log_LX3_thm}
Let $(X, \Delta)$ be a log Fano pair and $(\sX, \sL)/\A^1$ be a normal, ample 
test configuration for 
$(X, -(K_X+\Delta))$ with $(\sX, \Delta_{\sX}+\sX_0)$ log canonical. 
Then there exists 
a normal, ample test configuration $(\sX^{\AC}, \sL^{\AC})/\A^1$ for 
$(X, -(K_X+\Delta))$ with $(\sX^{\AC}, \Delta_{\sX^{\AC}}+\sX^{\AC}_0)$ 
log canonical such that: 
\begin{enumerate}
\renewcommand{\theenumi}{\arabic{enumi}}
\renewcommand{\labelenumi}{(\theenumi)}
\item\label{log_LX3_thm1}
$-(K_{\sX^{\AC}}+\Delta_{\sX^{\AC}})\sim_\Q\sL^{\AC}$. 
\item\label{log_LX3_thm2}
$\DF_\Delta(\sX^{\AC}, \sL^{\AC})\leq \DF_\Delta(\sX, \sL)$ holds. Moreover, 
equality holds if and only if $(\sX^{\AC}, \sL^{\AC})\simeq(\sX, \sL)$ over $\A^1$. 
\item\label{log_LX3_thm3}
$\Ding_\Delta(\sX^{\AC}, \sL^{\AC})\leq \Ding_\Delta(\sX, \sL)$ holds. Moreover, 
equality holds if and only if $(\sX^{\AC}, \sL^{\AC})\simeq(\sX, \sL)$ over $\A^1$. 
\item\label{log_LX3_thm4}
For any $\delta\in[0,1]$, we have the inequality 
\[
\Ding_\Delta(\sX, \sL)-\delta\cdot J^{\NA}(\sX, \sL)\geq
\Ding_\Delta(\sX^{\AC}, \sL^{\AC})-\delta\cdot J^{\NA}(\sX^{\AC}, \sL^{\AC}).
\]
\end{enumerate}
\end{thm}

\begin{proof}
This is proved in the same way as Proposition \ref{LX3_thm}. 
\end{proof}

\begin{thm}[{see Theorem \ref{LX4_thm}}]\label{log_LX4_thm}
Let $(X,\Delta)$ be a log Fano pair and 
$(\sX, -(K_{\sX/\A^1}+\Delta_{\sX}))/\A^1$ be a normal, ample 
test configuration for $(X, -(K_X+\Delta))$ with $(\sX, \Delta_{\sX}+\sX_0)$ log canonical. 
Then there exist $d\in\Z_{>0}$ and a projective birational 
$\G_m$-equivariant birational map $\sX^{(d)}
\dashrightarrow\sX^{\s}$ over $\A^1$ 
such that: 
\begin{enumerate}
\renewcommand{\theenumi}{\arabic{enumi}}
\renewcommand{\labelenumi}{(\theenumi)}
\item\label{log_LX4_thm1}
$(\sX^{\s}, -(K_{\sX^{\s}/\A^1}+\Delta_{\sX^{\s}}))/\A^1$ 
is a special test configuration for $((X, \Delta), -(K_X+\Delta))$. 
\item\label{log_LX4_thm2}
$\DF_\Delta(\sX^{\s}, -(K_{\sX^{\s}/\A^1}+\Delta_{\sX^{\s}}))\leq 
d\cdot\DF_\Delta(\sX, -(K_{\sX/\A^1}+\Delta_{\sX}))$ holds. 
Moreover,  
equality holds if and only if $(\sX, -(K_{\sX/\A^1}+\Delta_{\sX}))/\A^1$ 
is a special test configuration for 
$((X,\Delta), -(K_X+\Delta))$. We remark that 
$\DF_\Delta(\sX^{\s}, -(K_{\sX^{\s}/\A^1}+\Delta_{\sX^{\s}}))=
\Ding_\Delta(\sX^{\s}, 
-(K_{\sX^{\s}/\A^1}+\Delta_{\sX^{\s}}))$ and 
$\DF_\Delta(\sX, -(K_{\sX/\A^1}+\Delta_{\sX}))=
\Ding_\Delta(\sX, -(K_{\sX/\A^1}+\Delta_{\sX}))$ hold. 
\item\label{log_LX4_thm3}
For any $\delta\in[0,1]$, we have the inequality 
\begin{eqnarray*}
&&d(\Ding_\Delta(\sX, -(K_{\sX/\A^1}+\Delta_{\sX}))-
\delta\cdot J^{\NA}(\sX, -(K_{\sX/\A^1}+\Delta_{\sX})))\\
&\geq&\Ding_\Delta(\sX^{\s}, -(K_{\sX^{\s}/\A^1}+\Delta_{\sX^{\s}}))
-\delta\cdot J^{\NA}(\sX^{\s}, -(K_{\sX^{\s}/\A^1}+\Delta_{\sX^{\s}})).
\end{eqnarray*}
\end{enumerate}
\end{thm}

\begin{proof}
We follow the proof of Proposition \ref{LX4_thm}. 
For a logarithmic version of \cite[Theorem 6 (2)]{LX}, we use the following 
Theorem \ref{LX6_thm}.  
\end{proof}

\begin{thm}[{see \cite[Theorem 6 (2)]{LX}}]\label{LX6_thm}
Let $0\in C$ be a germ of a smooth curve and $C^*:=C\setminus\{0\}$. Let 
$f\colon\sX\to C$ be a projective morphism and $\Delta$ be a horizontal effective 
$\Q$-divisor on $\sX$ such that $-(K_{\sX}+\Delta)$ is ample over $C$, 
$(\sX, \Delta+\sX_0)$ is log canonical and $(\sX^*, \Delta^*)$ is klt, 
where $\sX^*:=f^{-1}(C^*)$ and $\Delta^*:=\Delta|_{\sX^*}$. Then there exists a 
finite morphism $\phi\colon (0\in C')\to (0\in C)$ between germs of smooth curves 
which is \'etale outside $0\in C$, a projective morphism $f^{\s}\colon \sX^{\s}\to C'$ 
and a horizontal effective $\Q$-divisor $\Delta^{\s}$ on $\sX^{\s}$ such that 
the following hold: 
\begin{itemize}
\item
$\sX^{\s}|_{(f^{\s})^{-1}(\phi^{-1}(C^*))}\simeq\sX^*\times_C C'$ over $\phi^{-1}(C^*)$. 
Moreover, $\Delta^{\s}$ is the closure of the pullback of $\Delta^*$. 
\item
$-(K_{\sX^{\s}}+\Delta^{\s})$ is ample over $C'$.
\item
The pair $(\sX^{\s}, \Delta^{\s}+\sX_0^{\s})$ is plt.
\item
The discrepancy $a(\sX_0^{\s}; \sX', \Delta')$ is equal to zero, 
where $(\sX', \Delta')$ is the normalization of the pullback of $(\sX, \Delta)$ to $\phi$.
\end{itemize}
\end{thm}

\begin{proof}
We give a proof for the readers' convenience. 
By the semistable reduction theorem, there exist a finite morphism 
$\phi\colon (0\in C')\to (0\in C)$ which is \'etale outside $0\in C$ 
(we set $\phi\colon\sX'\to\sX$ the normalization of $\sX\times_CC'$ and 
$\Delta':=\phi^*\Delta$)
and a projective birational morphism $\pi\colon\sY\to \sX'$ such that the pair 
$(\sY, \sY_0+\Exc(\pi)+\pi^{-1}_*\Delta')$ is a simple normal crossing pair and $\sY_0$ 
is reduced. 
We write 
\[
\pi^*(K_{\sX'{}^*}+\Delta'{}^*)+F^*=K_{\sY^*}+\pi_*^{-1}\Delta'{}^*+E^*,
\]
where $\sX'{}^*:=(f^{\s})^{-1}(\phi^{-1}(C^*))$, $\sY^*:=\pi^{-1}(\sX'{}^*)$, 
$\Delta'{}^*:=\Delta'|_{\sX'{}^*}$ and 
$E^*$, $F^*\geq 0$ without common components. Since 
$(\sX'{}^*, \Delta'{}^*)$ is klt, we have $\lceil E^*\rceil=0$. We can write 
\[
\pi^*(K_{\sX'}+\Delta')+F+B=K_{\sY}+\pi_*^{-1}\Delta'+E,
\]
where $E$, $F$ be the closures of $E^*$, $F^*$ and $B$ is vertical. 
Since $(\sX', \Delta'+\sX'_0)$ is log canonical, we have $B\geq 0$. 

Take a $\pi$-ample $\Q$-divisor $A$ on $\sY$ 
such that the support of $A$ is equal to the union 
of $\Exc(\pi)$ and the support of $\sY_0$. 
Let $G$ be the reduced divisor whose support is equal to the union of 
$\pi$-exceptional, horizontal divisors. Then 
$\delta G+\varepsilon A^h\geq 0$ and the support of $\delta G+\varepsilon A^h$ 
is equal to the support of $G$ for $0<\varepsilon\ll\delta\ll 1$, where $A^h$ is the 
horizontal part of $A$ and $A^v$ is the vertical part of $A$. 

Set $\sL_{\sY}:=\pi^*(-(K_{\sX'}+\Delta'))+\varepsilon A$. $\sL_{\sY}$ is ample 
over $C'$. Note that 
\[
K_{\sY}+\pi^{-1}_*\Delta'+E+\delta G+\sL_{\sY}\sim_{\Q, C'} (F+\delta G+
\varepsilon A^h)+(B+\varepsilon A^v). 
\]
Set $B+\varepsilon A^v=:\sum_{i=1}^ma_iE_i$. 
After a suitable small perturbation of the coefficients of $A$, we can assume that 
$a_1<\cdots<a_m$ and $a(E_1; \sX', \Delta')=0$. 
By \cite[Proposition 1]{LX}, we can run $(K_{\sY}+\pi^{-1}_*\Delta'+E+\delta G
+\sL_{\sY})$-MMP
\[
\sY=\sY^0\dashrightarrow\sY^1\dashrightarrow\cdots\dashrightarrow\sY^k
\]
over $C'$ with scaling $\sL_{\sY}$ such that $K_{\sY^k}+\Delta^k+\sL^k\sim_{\Q, C'}
0$, $\sL^k$ is nef and big over $C'$ and the support of the exceptional divisors of 
the rational map $\sY\dashrightarrow\sY^k$ is equal to the union of the support of 
$G$ and $\bigcup_{i=2}^mE_i$, where 
$\Delta^k$, $\sL^k$ is the strict transform of $\pi^{-1}_*\Delta'$, $\sL_{\sY}$ 
on $\sY^k$, respectively. 
From the base-point freeness, $\sL^k$ is semiample 
over $C'$. Let us consider the canonical model $\psi\colon\sY^k\to\sX^{\s}$ of $\sL^k$
over $C'$. Obviously, $\sX^{\s}$ coincides with $\sX'$ outside $0\in C'$. Moreover, 
$\sX_0^{\s}$ is the image of $E_1$. Since the pair $(\sY, \pi^{-1}_*\Delta'+E+\delta G
+\sY_0)$ is dlt and $\sY_0^k$ is irreducible, 
the pair $(\sY^k, \Delta^k+\sY_0^k)$ is plt. Thus the pair 
$(\sX^{\s}, \Delta^{\s}+\sX_0^{\s})$ is also plt. 
\end{proof}

The following is a direct consequence of Theorems \ref{log_LX2_thm}, 
\ref{log_LX3_thm} and \ref{log_LX4_thm}. See also \cite[Theorem 7]{LX}, 
\cite[Theorem 2.1]{BBJ} and \cite[Theorem 8.6]{BHJ}. 

\begin{corollary}[{see Corollary \ref{stc_cor}}]\label{log_stc_cor}
Let $(X, \Delta)$ be a log Fano pair. 
\begin{enumerate}
\renewcommand{\theenumi}{\arabic{enumi}}
\renewcommand{\labelenumi}{(\theenumi)}
\item\label{log_stc_cor1}
For any $\delta\in[0,1)$, the following conditions are equivalent: 
\begin{enumerate}
\renewcommand{\theenumii}{\roman{enumii}}
\renewcommand{\labelenumii}{(\theenumii)}
\item\label{log_stc_cor11}
$\DF_\Delta(\sX, \sL)\geq \delta\cdot J^{\NA}(\sX, \sL)$ holds for any normal, ample 
test configuration $(\sX, \sL)/\A^1$ for $(X, -(K_X+\Delta))$. 
\item\label{log_stc_cor12}
$\Ding_\Delta(\sX, \sL)\geq \delta\cdot J^{\NA}(\sX, \sL)$ holds for any normal, ample 
test configuration $(\sX, \sL)/\A^1$ for $(X, -(K_X+\Delta))$. 
\item\label{log_stc_cor13}
$\DF_\Delta(\sX, -(K_{\sX/\A^1}+\Delta_{\sX}))\geq 
\delta\cdot J^{\NA}(\sX, -(K_{\sX/\A^1}+\Delta_{\sX}))$ holds\, 
for\, any\, special\, test\, configuration\, $(\sX, -(K_{\sX/\A^1}+\Delta_{\sX}))/\A^1$ 
for $((X,\Delta), -(K_X+\Delta))$. 
\end{enumerate}
\item\label{log_stc_cor2}
$(X, \Delta)$ is K-stable $($resp.\ K-polystable$)$ if and only if 
$(X, \Delta)$ is Ding stable $($resp.\ Ding polystable$)$. 
\end{enumerate}
\end{corollary}

The following theorem is a logarithmic version of Theorem \ref{Kv_thm}. 
See also \cite[Propositions 4.3 and 4.5]{LL}. 

\begin{thm}[{see Theorem \ref{Kv_thm}}]\label{log_Kv_thm}
Let $(X, \Delta)$ be a log Fano pair. Assume that there exists 
$\delta\in[0,1)$ such that 
$\Ding_\Delta(\sX, \sL)\geq \delta \cdot J^{\NA}(\sX, \sL)$ holds for any 
normal, semiample test configuration $(\sX, \sL)/\A^1$ for \,
$(X, -(K_X+\Delta))$. Then 
$\beta_{(X, \Delta)}(F)\geq\delta\cdot j_{(X, \Delta)}(F)$ holds 
for any prime divisor $F$ over $X$. 
\end{thm}

\begin{proof}
This is proved in the same way as Proposition \ref{Kv_thm}. 
\end{proof}

The following theorem is a logarithmic version of Theorem \ref{vK_thm}. 
See also \cite[Lemma 4.8]{LX2}.

\begin{thm}[{see Theorem \ref{vK_thm}}]\label{log_vK_thm}
Let $(X, \Delta)$ be a log Fano pair and let $(\sX, -(K_{\sX/\A^1}+\Delta_{\sX}))/\A^1$ 
be a normal, ample test configuration for $(X, -(K_X+\Delta))$ with $\sX_0$ irreducible 
and reduced. Then the divisorial valuation $v_{\sX_0}$ on $X$ defined in 
Proposition \ref{discrep_prop} is dreamy over $(X, \Delta)$ 
and we have the equalities 
\[
\DF_\Delta(\sX, -(K_{\sX/\A^1}+\Delta_{\sX}))
=\frac{\beta_{(X, \Delta)}(v_{\sX_0})}{\vol_X(-(K_X+\Delta))}
\] 
and 
\[
J^{\NA}(\sX, -(K_{\sX/\A^1}+\Delta_{\sX}))
=\frac{j_{(X, \Delta)}(v_{\sX_0})}{\vol_X(-(K_X+\Delta))}.
\] 
\end{thm}

\begin{proof}
This is proved in the same way as Proposition \ref{vK_thm}. 
\end{proof}

\begin{proof}[Proof of Theorem \ref{log_mainthm}]
This is an immediate consequence of Corollary \ref{log_stc_cor}, 
Theorems \ref{log_Kv_thm} and \ref{log_vK_thm}.
\end{proof}

Finally, we interpret Corollary \ref{log_stc_cor} in the non-Archimedean 
language introduced in \cite{BHJ}. 
From now on, let us fix a log Fano pair $(X, \Delta)$ of dimension $n$. 

\begin{definition}[{\cite[Definitions 6.1 and 6.2]{BHJ}}]\label{NA_dfn}
\begin{enumerate}
\renewcommand{\theenumi}{\arabic{enumi}}
\renewcommand{\labelenumi}{(\theenumi)}
\item\label{NA_dfn1}
Two \,\,test configurations $(\sX^1, \sL^1)/\A^1$ and $(\sX^2, \sL^2)/\A^1$ for 
$(X, -(K_X+\Delta))$ are said to be \emph{equivalent} if there is 
$(\sX^3, \sL^3)/\A^1$ that is a pullback of $(\sX^1, \sL^1)/\A^1$ and 
$(\sX^2, \sL^2)/\A^1$. 
\item\label{NA_dfn2}
A \emph{non-Archimedean metric} on $-(K_X+\Delta)$ is an equivalence class of a test 
configurations for $(X, -(K_X+\Delta))$. The equivalence class of $(X\times\A^1, p_1^*
(-(K_X+\Delta)))/\A^1$ is denoted by $\phi_{\triv}$. 
\item\label{NA_dfn3}
A non-Archimedean metric $\phi$ on $-(K_X+\Delta)$ is said to be \emph{positive} 
if some $(\sX, \sL)/\A^1$ which gives $\phi$ is a semiample test configuration. 
The set of non-Archimedean positive metrics on $-(K_X+\Delta)$ is denoted by 
$\sH^{\NA}(-(K_X+\Delta))$.
\end{enumerate}
\end{definition}

\begin{definition}[{see \cite[\S 6.7]{BHJ}}]\label{MA_dfn}
Let $X^{\DIV}$ be the set of divisorial valuations on $X$. For any 
$\phi\in\sH^{\NA}(-(K_X+\Delta))$, we obtain a finite atomic measure 
$\MA^{\NA}(\phi)$, named a \emph{Monge-Amp\`ere measure}, on $X^{\DIV}$ as 
follows. Take any representation $(\sX, \sL)/\A^1$ of $\phi$ with $\sX$ normal, 
and set 
\[
\MA^{\NA}(\phi):=\frac{1}{((-(K_X+\Delta))^{\cdot n})}
\sum_{\substack{E\subset\sX_0\\ \text{irreducible component}}}
m_E\left(\sL|_E^{\cdot n}\right)\delta_{v_E}, 
\]
where $m_E:=\ord_E(\sX_0)$ and $v_E\in X^{\DIV}$ be as in \S 
\ref{val_section}.
\end{definition}

Consider the non-Archimedean positive metric $\phi\in\sH^{\NA}(-(K_X+\Delta))$ 
associated to a \emph{special} test configuration for $((X, \Delta), -(K_X+\Delta))$. 
Then there exists a dreamy prime divisor $F$ over $(X, \Delta)$ such that 
$\MA^{\NA}(\phi)$ is a multiple of $\delta_{\ord_F}$. The invariant 
$\beta_{(X, \Delta)}(F)$ (resp.\ $j_{(X, \Delta)}(F)$) is equal to 
$((-(K_X+\Delta))^{\cdot n})$ times 
the non-Archimedean 
Ding functional (resp.\ $J$-functional) evaluated at $\phi$ (see \cite[\S 7]{BHJ}). 
Therefore we get the following: 

\begin{corollary}\label{NA_cor}
It suffies to consider non-Archimedean metrics whose 
Monge-Amp\`ere measure is a Dirac mass at the divisorial valuation associated to 
dreamy prime divisors over $(X, \Delta)$ in order to test uniform K-stability $($resp.\ K-semistability$)$ of $(X, \Delta)$.
\end{corollary}

\end{document}